\documentclass[11pt]{article}
\usepackage{graphicx}
\usepackage{amsmath,amssymb}
\usepackage{amsthm}
\usepackage{epstopdf}
\usepackage{enumerate}
\usepackage{cases}
\usepackage{mathrsfs}
\usepackage[margin=1.0cm]{geometry}
\usepackage{color}
\definecolor{color1}{rgb}{0.7,0.4,0.4}
\definecolor{color2}{rgb}{0.4,0.4,0.7}
\pagecolor{white}
\color{black}

\numberwithin{equation}{section}

\usepackage{hyperref}
\hypersetup{
	colorlinks=true,
	urlcolor =blue,
	linkcolor=color2,
	citecolor=color1
}

\usepackage{appendix}
\usepackage{setspace}
\newcommand{\prob}[1]{\mathbb{P}\left[#1\right]}
\newcommand{\expected}[1]{\mathbb{E}\left[#1\right]}
\newcommand{\expectedcond}[2]{\mathbb{E}^{(#2)}\left[#1\right]}

\newcommand{\sgn}[1]{\operatorname{sgn}\left(#1\right)}
\renewcommand{\bar}[1]{\overline{#1}}
\newcommand{\dd}{\partial}

\title{\bf \sc \Huge Stochastic averaging of dynamical systems with multiple time scales forced with $\alpha$-stable noise}
\author{William F. Thompson\footnote{Dept. of Mathematics, U. of British Columbia, Vancouver, BC, Canada. Corresponding author: [ \href{mailto:wft@math.ubc.ca}{\nolinkurl{wft@math.ubc.ca}} ]} , Rachel A. Kuske \footnote{Dept. of Mathematics, U. of British Columbia, Vancouver, BC, Canada.} , Adam H. Monahan\footnote{School of Earth and Ocean Science, U. of Victoria, Victoria, BC, Canada.}}
\date{}

\begin{document}
\maketitle

\begin{abstract}
Stochastic averaging allows for the reduction of the dimension and complexity of stochastic dynamical systems with multiple time scales, replacing fast variables with statistically equivalent stochastic processes in order to analyze variables evolving on the slow time scale. These procedures have been studied extensively for systems driven by  Gaussian noise, but very little attention has been given to the case of $\alpha$-stable noise forcing which arises naturally from heavy-tailed stochastic perturbations. In this paper, we study nonlinear fast-slow stochastic dynamical systems in which the fast variables are driven by additive $\alpha$-stable noise perturbations, and the slow variables depend linearly on the fast variables. Using a combination of perturbation methods and Fourier analysis, we derive stochastic averaging approximations for the statistical contributions of the fast variables to the slow variables. In the case that the diffusion term of the reduced model depends on the state of the slow variable, we show that this term is interpreted in terms of the Marcus calculus. For the case $\alpha = 2$, which corresponds to Gaussian noise, the results are consistent with previous results for stochastic averaging in the Gaussian case. Although the main results are derived analytically for $1 < \alpha < 2$, we provide evidence of their validity for $\alpha < 1$ with numerical examples. We numerically simulate both linear and nonlinear systems and the corresponding reduced models demonstrating good agreement for their stationary distributions and temporal dependence properties.
\end{abstract}

\section{Introduction}
\label{sec:intro}

Stochastic averaging (also referred to as ``stochastic reduction'' or in some cases, ``stochastic homogenization'') refers to techniques whereby a stochastic dynamical system with multiple time scales is approximated by a reduced dynamical system on a subset of slow variables sharing essential statistical properties of the slow variables of the original system. These methods allow the modelling of the behaviour in the weak sense (i.e. distributions, moments) of processes influenced by faster processes without having to explicitly model the dynamical details of the faster variables. This averaging can reduce computational load and the resulting approximation is typically easier to analyze and interpret. Stochastic averaging methods have been discussed widely and we mention a few relevant references here. Freidlin and Wentzell \cite{Freidlin1984} provide a mathematically rigorous overview of fundamental stochastic averaging procedures, such as those developed by Khas'minskii \cite{Khasminskii1966a} and extended by Papanicolaou  \cite{Papanicolaou1974} and Borodin \cite{Borodin1977} . Majda, Timofeyev and Vanden-Ejinden developed a generalized stochastic averaging procedure for quadratically nonlinear dynamical systems with slow resolved components and fast unresolved components (MTV theory) \cite{ Majda2001,Majda2002, Majda2003}. Work by Givon et al. \cite{Givon2004} and Pavilotis \cite{Pavliotis2007} cover stochastic averaging and homogenization results obtained by perturbation analysis of the Backward Kolmogorov equation in the context of Gaussian noise.  

MTV theory was developed primarily in the context of weather-climate interaction which is a key application for stochastic averaging theory. Variability in the atmosphere-ocean-ice-land surface system occurs over a range of time scales spanning many orders of magnitude \cite{Palmer2008}. Generally speaking, the slowly-evolving climate variables describe the state of the ocean, the land surface, ice, and atmospheric composition, while the fast weather variables describe the circulation and thermodynamic state of the atmosphere \cite{Saltzman2002}. In order to model longer-term properties of the climate system, evolving on time-scales of years to centuries, it is undesirable to explicitly simulate meteorological processes (i.e. weather) which evolve on much faster time-scales ranging from hours to days to weeks. However, variability of the climate and weather systems are mutually dependent.

There are three established stochastic averaging approximations that are generalized in this paper. We refer to them using the terminology introduced in \cite{Arnold2003, Mona11} so as to be consistent with previous literature. Given a stochastic dynamical system with Gaussian white noise forcing and variables that evolve on both fast and slow time scales, the (A) approximation for the mean of the slow dynamics of the system is determined by replacing any functions of the fast variables with their means conditioned on the value(s) of the slow variable(s). This approximation accounts for the influence of functions of the fast variables only through their means, and otherwise neglects their variability. The (L) approximation improves upon the (A) approximation by adding a correction in the form of an Ornstein-Uhlenbeck-like process with dynamics determined via an analysis of the full system about the mean of the slow variables. As (L) describes local variability around (A), it is unable to account for transitions between metastable states. Both the (A) and (L) approximations are derived in the work of Khas'minskii \cite{Khasminskii1966b, Khasminskii1966a}. Lastly, we consider the (N+) approximation which represents the slow dynamics as a more general stochastic differential equation (SDE) with potentially nonlinear and/or with state-dependent noise. The (N+) approximation was inspired by the work of Hasselmann \cite{Hasselmann1976}, with more rigorous analysis of the approximation in \cite{Kifer2003}. This approximation was derived formally in \cite{Hasselmann1976} by considering an arbitrary point in the state space of the slow variables and performing a local analysis of the dispersion at that point to derive drift and diffusion coefficients. Other studies offer derivations of the (N+) approximation via projection methods and perturbation analysis of the Forward Kolmogorov equation \cite{Just2003,Just2001}.

Previous derivations of the (L) and (N+) approximations rely on the evaluation of the integrated autocovariance function of the fast process conditioned on the state of the slow variable in the presence of Gaussian noise forcing. When the autocovariance function is not defined, if for example the second moment of the stochastic forcing does not exist, the (L) and (N+) approximations cannot be defined in the usual way as given in \cite{Arnold2003, Mona11}. This paper focuses on averaging results for systems driven by processes with increments that are distributed according to an $\alpha$-stable law, for which the autocovariance is not defined \cite{Feller1966II} (except for the case when $\alpha = 2$). A stochastic averaging procedure for systems with multiple time scales and $\alpha$-stable noise forcing is relevant for modelling problems with heavy-tailed statistics, as a shifted and normalized series of heavy-tailed random variables converges in distribution to an $\alpha$-stable random variable \cite{Taqqu1994}. Hence, they form a family of attracting distributions for independent, identically-distributed (i.i.d.) sums of heavy-tailed random variables, analogous to the Gaussian distribution as the limiting distribution for i.i.d. sums of random variables with finite variance. An integrated $\alpha$-stable noise process is characterized by its ``superdiffusive'' behaviour whereby the observed variance in the trajectories grows superlinearly in time and are sometimes referred to as a ``L\'{e}vy flights'' in the physics literature. The dynamics of atmospheric tracers have been shown to be superdiffusive in several contexts \cite{Huber2001,Seo2000}. Similar dynamics have also been predicted to be generated by chaotic dispersion and shear flows \cite{Solomon1994}. Processes such as these would be considered relatively fast atmospheric processes. In addition to these results, an analysis by Ditlevsen has suggested the presence of $\alpha$-stable forcing on longer time scales is present in time series of calcium concentration from Greenland ice cores \cite{Ditlevsen1999, Hein2008} recording high latitude climate variability during the most recent glacial period. 

An $\alpha$-stable noise process is one for which the increments are $\alpha$-stable random variables whose distribution we denote by $\mathcal{S}_{\alpha}(\beta,\sigma)$. The distribution depends on three parameters: the stability index $\alpha \in (0,2]$ which determines the rate of decay of the tails of the distribution, the skewness parameter $\beta \in [-1,1]$ which describes the asymmetry in the distribution, and the scale parameter $\sigma > 0$ which governs the ``width'' of the density function \cite{Applebaum2004, Taqqu1994}. The density function for an $\alpha$-stable random variable cannot be expressed in closed form, except for a small number of special cases. Instead, an $\alpha$-stable random variable, $L \sim \mathcal{S}_{\alpha}(\beta,\sigma)$, is typically characterized through its characteristic function (CF), $\Phi_{L}(k)$,
\begin{equation}
\Phi_{L}(k) = \exp\left(-\sigma^{\alpha}|k|^{\alpha}\Xi(k;\alpha,\beta)\right) \label{alpha_stable_characteristic function}
\end{equation}
where
\begin{equation}
\Xi(s;\alpha,\beta) = 1 + i\beta\sgn{s}\tan(\pi\alpha/2). \label{Xi}
\end{equation}
It is easy to see that $\alpha = 2$ in (\ref{alpha_stable_characteristic function}) corresponds to the special case of a Gaussian random variable with variance $2\sigma^{2}$ and mean $0$. The parameter $\beta$ is inconsequential in this case, since $\Xi(k;2,\beta) = 1$ for all $k,\beta$. If $\alpha < 2$, the second moment of the distribution does not exist. When $\alpha < 1$, the first moment does not exist. For $\alpha$-stable increments in continuous time denoted by $dL_{t}^{(\alpha,\beta)}$, the scale parameter of the increment is nonlinearly related to that of the infinitesimal timestep $dt$. Formally this relationship implies that $\sigma^{\alpha} \propto dt$ and so the noisy increments $dL_{t}^{(\alpha,\beta)} \sim \mathcal{S}_{\alpha}\left(\beta,(dt)^{1/\alpha}\right)$. If $\beta = 0$, the noise is symmetric. Otherwise, it is asymmetric.

The mathematical problem of stochastic averaging of heavy-tailed processes for linear systems has been investigated in \cite{Srokowski2011} in the context of the Forward Kolmogorov equation (FKE). We use the more general FKE terminology throughout the paper, noting that if the stochastic forcing is Gaussian, the equation is usually called the Fokker-Planck equation. More recent work by Chechkin and Pavlyukevich \cite{Chechkin2014} and Li et al. \cite{Li2014} give Wong-Zakai type results for systems driven by a fast compound Poisson jump processes. Results like those of Xu et al. \cite{Xu2011} consider strong-sense (i.e pathwise) stochastic averaging of $\alpha$-stable forced systems, determining reduced SDE models that approximate the full systems in a pathwise sense. They consider a class of two-time scale systems in which the fast dynamics corresponds to an explicit time-dependence (on fast timescales) of the slow variable's drift and diffusion. We consider a much broader class of systems, but only in the weaker sense of convergence in distribution.

\begin{figure}[h]
\centering
\includegraphics[width=\textwidth]{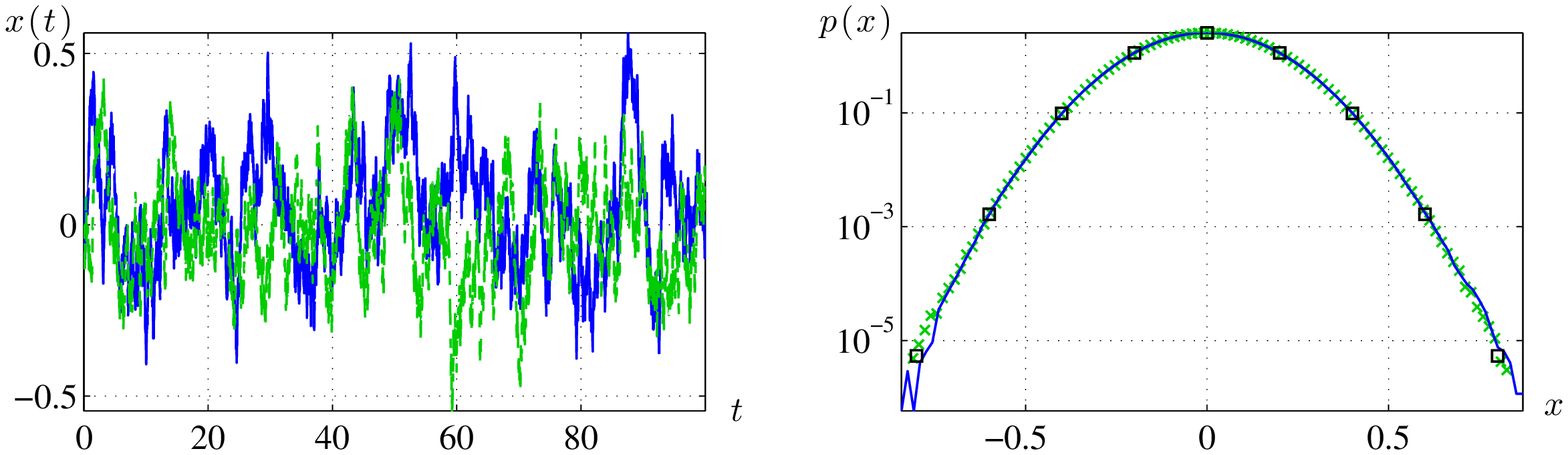}\\
\includegraphics[width=\textwidth]{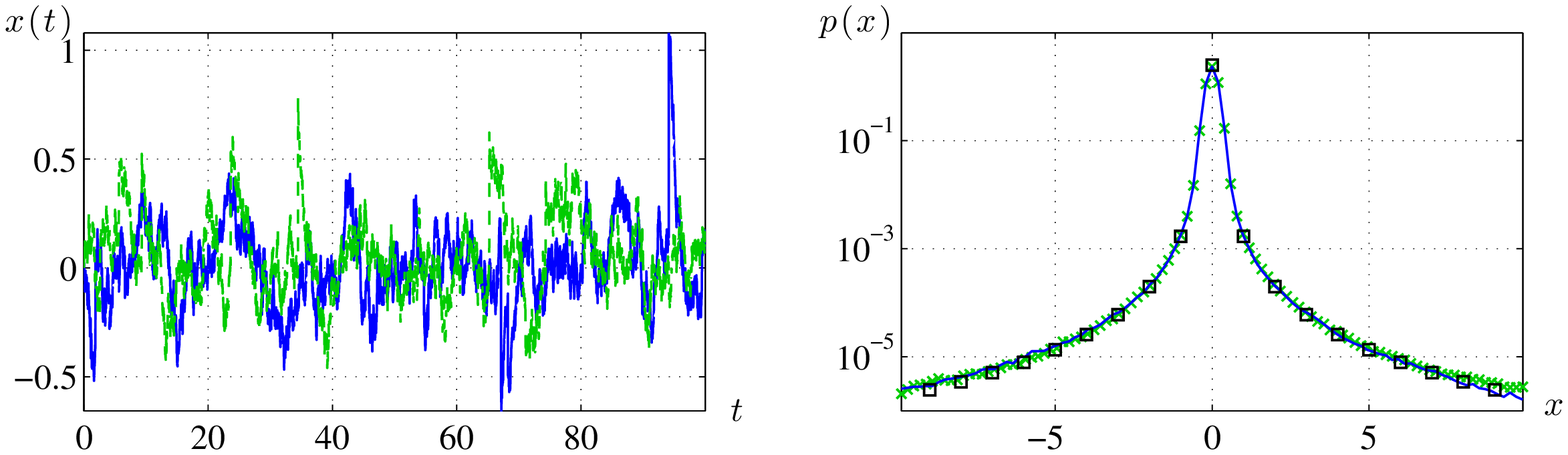}
\caption{A comparison of the time series and simulated marginal distributions of the fast-slow system (\ref{linsys1},  \ref{linsys2}) and the stochastically averaged system (\ref{Lapprox_linear}). Top-Left: A sample of the time series for the slow variable of both systems: $x_{t}$ (solid blue) and $\bar{x}_{t} + \xi_{t}$ (dashed green). Top-Right: Normalized histogram-based estimates for the density for the slow variable. Blue line: full system. Green crosses: (L) approximation. Black squares: Numerically evaluated probability density function for $\xi$ from the (L) approximation (\ref{Lapprox_linear}) using the MATLAB package STBL as described in \S \ref{sec:comp_results}. Parameters for this simulation are $(a = 0.2,\,b = 0.7,\,c = 1,\,\epsilon = 0.01,\,\alpha = 2.0,\,x_{0} = 0)$. The bottom panels are as for the top, but with $\alpha$-stable noise where $\alpha = 1.9$, $\beta = 0$.} \label{linhisto_intro}
\end{figure}

We derive stochastic averaging approximations in the weak (distributional) sense for fast-slow stochastic dynamical systems with additive $\alpha$-stable noise forcing, obtaining reduced dynamical approximations for the slow component of the full fast-slow system. As in \cite{Srokowski2011}, we analyze the FKE of the fast-slow systems, but in contrast we consider both linear and nonlinear systems under appropriate transformations to allow an analysis of the FKE in Fourier space. This approach allows us to consider a broader class of systems than has been investigated previously in the context of stochastic averaging of $\alpha$-stable noise processes. These systems are linear in the fast variable, while allowing for some nonlinear interactions between slow and fast variables and fully nonlinear behaviour in the dynamics of the slow variables. The averaging approximations are conceptually similar to existing stochastic averaging methods used for systems with Gaussian white noise \cite{Arnold2003, Mona11}, but the approach used for the $\alpha$-stable case is significantly different from that used for the Gaussian case to determine the appropriate stochastic forcing for the reduced dynamics. The analysis here is based on the joint CF of the fast component and a fluctuation about the mean of the slow variable. For the (N+) approximations, Marcus' stochastic calculus \cite{Chechkin2014, Marcus1978} for jump processes is required if the approximation has a state-dependent ``diffusion'' term. Analogous to Stratonovich calculus for Gaussian processes, the Marcus interpretation of stochastic integrals allows for the rules of ordinary calculus, like the chain rule, to hold. In Fig. \ref{linhisto_intro}, we show an example of stochastic averaging of a linear system where $\alpha = 2$ and $\alpha < 2$. For such systems, the (L) and (N+) approximations are the same and they compare well with the full system.

In \S \ref{sec:averaging}, we define the stochastic averaging problem and write the expression for the (A) approximation in the case where the stochastic forcing is $\alpha$-stable noise, rather than Gaussian. In \S \ref{subsec:Lapprox} and \S \ref{subsec:N+approx}, we derive the analogous (L) and (N+) approximations for systems satisfying the appropriate conditions, demonstrating that the noise terms in the (N+) approximation are interpreted in the sense of Marcus. Since the Marcus calculus is less well-known, we include a brief discussion in Appendix \ref{app:simulation_Marcus}. In \S \ref{sec:comp_results}, we apply our approximations to one linear system and three distinct nonlinear systems, comparing the averaging approximations to the full systems. These numerical results show agreement between the approximations and the full systems not only in probability distribution, but also in the autocodifference which characterizes the memory of the process. Our conclusions are outlined in \S \ref{sec:conc}.

\section{Stochastic averaging and $\alpha$-stable noise}
\label{sec:averaging}

We study the canonical problem of stochastic averaging of a fast-slow system with the stochastic forcing given by the increments of an $\alpha$-stable L\'{e}vy process, $L_{t}^{(\alpha,\beta)}$,
\begin{align}
dx_{t} &= f(x_{t},y_{t})\,dt, \label{canon1}\\
\epsilon\,dy_{t} &= g(x_{t},y_{t})\,dt + \epsilon^{\gamma}b\,dL_{t}^{(\alpha,\beta)},\label{canon2}
\end{align}
where $\gamma$ and $b$ are constants and $x_{0},\,y_{0}$ are known. The parameter $\epsilon$ is a small positive number ($0 < \epsilon \ll 1$) that characterizes the well-separated time scales between the slow variable $x_{t}$ and the fast variable $y_{t}$. The variable $y_{t}$ is a stationary and ergodic stochastic process forced with an additive $\alpha$-stable L\'{e}vy noise $dL_{t}^{(\alpha,\beta)}$ with increments distributed as $\mathcal{S}_{\alpha}\left( \beta,(dt)^{1/\alpha}\right)$. In this study, we restrict our analysis to the case $\alpha \in (1,2]$ with $\expected{dL_{t}^{(\alpha,\beta)}} = 0$. The system (\ref{canon1}, \ref{canon2}) includes the well-studied case of Gaussian noise forcing, corresponding to $\alpha = 2$. When $\alpha = 2$, we replace $dL_{t}^{(2,\beta)}$ with $\sqrt{2}\,dW_{t}$, where $W_{t}$ is a standard Wiener process and the factor of $\sqrt{2}$ follows from the definition of the $\alpha$-stable CF (\ref{alpha_stable_characteristic function}). Since $\alpha$-stable noise processes are best analyzed through their CFs, we consider forms of $f,g$ that are linear in $y_{t}$,
\begin{equation}
f(x,y) = f_{1}(x) + \epsilon^{-\gamma}f_{2}(x)y, \quad g(x,y) = \epsilon^{\gamma}g_{1}(x) + g_{2}(x)y, \label{drifts}
\end{equation}
to make analysis in Fourier space straightforward. Specifically, the linear form allows us to give results in terms of closed form functions of the parameters. Some preliminary work for examples that include polynomials of $y$ indicates the limited availability of explicit analytical expressions. Within the form (\ref{drifts}), we also assume $g_{2}(x) < 0$ for all $x$ to ensure contracting dynamics so that the distribution of $y|x$ reaches conditional stationary behaviour on a fast time scale. For the (N+) approximation in \S \ref{subsec:N+approx}, we require also that $f_{2}(x)$ has the same sign for all $x$ and is never equal to zero in order to define a particular monotonic transformation.

The goal of stochastic averaging is to derive an approximation in the weak sense to the SDE of $x_{t}$ based on the distributional properties of $y_{t}$. We obtain three different stochastic averaging approximations resulting in averaged deterministic, linear stochastic, and nonlinear stochastic models, as described below.

\subsection{(A) approximation}
The (A) approximation is an approximation of $x_{t}$ in the mean and is therefore deterministic. It is given by $\bar{x}_{t}$ where
\begin{equation}
d\bar{x}_{t} = \bar{f}(\bar{x}_{t})\,dt, \quad \bar{x}_{0} = x_{0} \label{Aapprox}
\end{equation}
and
\begin{equation}
\bar{f}(x) = \expectedcond{f(x,y)}{y|x} = \int_{-\infty}^{\infty}f(x,y)p(y|x)\,dy. \label{fbar}
\end{equation}
The notation $\expectedcond{f(x,y)}{y|x}$ indicates the expected value of $f(x,y)$ conditioned on a fixed $x$ and the term $p(y|x)$ is the conditional density for $y_{t}$ for a fixed $x_{t} = x$. For this approximation we have implicitly assumed that $\bar{f}$ can be evaluated for $x_{t}, y_{t}$ in (\ref{canon1}, \ref{canon2}). Indeed, by our assumption that $\alpha \in (1,2]$ and the form of $f,g$ given in (\ref{drifts}) holds, then $\bar{f}$ is evaluated using (\ref{fbar}) giving us $\expectedcond{y}{y|x} = -\epsilon^{\gamma}g_{1}(x)/g_{2}(x)$ and thus,
\begin{equation}
\bar{f}(z) = f_{1}(z) - \frac{f_{2}(z)g_{1}(z)}{g_{2}(z)}. \label{fbar_specific}
\end{equation}
For the Gaussian noise case, the (A) approximation has been studied in detail \cite{Arnold2003, Freidlin1984,Khasminskii1966a}. In the next section, we show that the (A) approximation holds in the case of $\alpha$-stable noise forcing if $\alpha > 1$ and $\gamma > 1 -1/\alpha$ in the strict $\epsilon \rightarrow 0$ limit. For $0 < \epsilon \ll 1$, a more accurate approximation called the (L) approximation emerges that includes corrections to the (A) approximation.

\subsection{(L) approximation}
\label{subsec:Lapprox}
The (L) approximation is an Ornstein-Uhlenbeck-L\'{e}vy process (OULp) obtained by linearizing about the mean $\bar{x}$ as defined in (\ref{Aapprox}). We first decompose the dynamics into a mean component and a fluctuation as in \cite{Mona11}. The properties of $\alpha$-stable noise require analysis of CFs to determine the appropriate stochastic approximation for the fluctuation.

We define the fluctuation about the mean as $\xi_{t} = (x_{t} - \bar{x}_{t})$ and derive its dynamics by substituting in (\ref{canon1}) and using (\ref{Aapprox}),
\begin{align}
d\xi_{t} &= [f(x_{t},y_{t}) - \bar{f}(\bar{x}_{t})]\,dt \\
&= \left[\bar{f}(x_{t}) - \bar{f}(\bar{x}_{t})\right]\,dt + \hat{f}(x_{t},y_{t})\,dt, \label{eqn:xi_dyn_exact}
\end{align}
where $\hat{f}(z,y) = f(z,y) - \bar{f}(z)$. Substituting $x_{t} = \bar{x}_{t} + \xi_{t}$ and keeping leading order contributions for $|\xi_{t}| \ll 1$ in (\ref{eqn:xi_dyn_exact}) yields the linear equation for $\xi_{t}$,
\begin{equation}
d\xi_{t} \simeq  \bar{f}'(\bar{x}_{t})\xi_{t}\,dt + \hat{f}(\bar{x}_{t},y_{t})\,dt. \label{lin_approx_fbar}
\end{equation}
To complete the (L) approximation, we find a stochastic approximation for $\hat{f}\,dt$ that has its distributional properties on the $t$ time scale. The variable $v_{t}$ captures the stochastic forcing in the $\xi$ equation. Its dynamics are
\begin{align}
dv_{t} &= \hat{f}(\bar{x}_{t},y_{t})\,dt, \quad v_{0} = 0, \quad \mbox{where} \quad \hat{f}(z,y) = f_{2}(z)\left(\epsilon^{-\gamma}y + \frac{g_{1}(z)}{g_{2}(z)}\right), \label{v_dynamics} 
\end{align}
and we have substituted (\ref{drifts}) and (\ref{fbar_specific}) into (\ref{lin_approx_fbar}) to obtain an explicit expression for $\hat{f}$. Rewriting (\ref{lin_approx_fbar}) using (\ref{v_dynamics}), we see the forcing effect of the increments $dv_{t}$ on the dynamics of $\xi$,
\begin{equation}
d\xi_{t} = \bar{f}'(\bar{x}_{t})\xi_{t}\,dt + dv_{t}, \quad \xi(0) = 0. \label{xi_dynamics_w_v}
\end{equation}
In the case $\alpha = 2$, the autocovariance $C(s;x)$ can be evaluated and the fluctuation $dv_{t}$ in (\ref{xi_dynamics_w_v}) is replaced with a Gaussian white noise forcing $\sigma(\bar{x})\,dW_{t}$ \cite{Arnold2003,Mona11}, where
\begin{equation}
\sigma^{2}(x) = \int_{-\infty}^{\infty}C(s;x)\,ds, \quad C(s;x) = \lim_{t \rightarrow \infty} \expectedcond{\hat{f}(x,y_{t+s})\hat{f}(x,y_{t})}{y|x}. \label{L_autoCov}
\end{equation}
This result is justified for $C(s;x)$ having a $\delta$-function behaviour similar to that of Gaussian white noise in the limit $\epsilon \rightarrow 0$ \cite{Mona11}. However, the autocovariance function $C(s;x)$ of $\hat{f}(x,y)$ is undefined when $\alpha < 2$. While moments with order greater than or equal to $\alpha$ do not exist for $\alpha$-stable distributions, the CF is defined for any $\alpha < 2$. Thus, we analyze the process $v_{t}$ in terms of its CF in order to determine an appropriate stochastic forcing approximation for $dv_{t}$ when $\alpha < 2$.

The joint probability density function (PDF) $p(y,v,t)$ for the system (\ref{canon2}), (\ref{v_dynamics}) satisfies the fractional FKE (also known as the Fokker-Planck equation if $\alpha = 2$),
\begin{equation}
\begin{cases} & \displaystyle \frac{\dd p}{\dd t} = - \frac{\dd}{\dd v}\left( \hat{f}(\bar{x},y)p\right) - \frac{1}{\epsilon}\frac{\dd}{\dd y}\left(g(\bar{x},y)p \right) + \frac{b^{\alpha}}{\epsilon^{(1-\gamma)\alpha}}\mathcal{D}_{y}^{(\alpha,\beta)}p, \\ &p(y,v,0) = \delta(v)\delta(y - y_{0}). 
\end{cases} \label{FPE_L_canon_reduced}
\end{equation}
where $\mathcal{D}_{y}^{(\alpha,\beta)}$ is the the Riesz-Feller derivative, a generalization of the Liouville-Riemann fractional derivative \cite{Chaves1998}, acting with respect to $y$ associated with the $\alpha$-stable noise process. This operator can be defined in terms of a Fourier transform, $F[u] = \int_{\mathbb{R}}\exp(ily)u(y)\,dy$, and its inverse:
\begin{align}
\mathcal{D}_{y}^{(\alpha,\beta)}u &= F^{-1}\left[\frac{-1}{\cos(\pi\alpha/2)}\left(\frac{1 + \beta}{2}(il)^{\alpha} + \frac{1 - \beta}{2}(-il)^{\alpha}\right)F\left[u\right]\right] \label{eqn:RieszFellerDeriv}
\end{align}
For $0 < \epsilon \ll 1$, we initially treat  $\bar{x}$ as a constant relative to the fast variables $v$ and $y$. Let $\psi$ be the CF for the $(y,v)$ system,
\begin{equation}
\psi(l,m,t) = \mathcal{F}[p] = \iint_{\mathbb{R}^{2}}\exp(ily + imv)p(y,v,t)\,dydv. \label{char_func}
\end{equation}
Taking the Fourier transform (FT) of (\ref{FPE_L_canon_reduced}), we note that the noise terms are entering the dynamics through $y$ and not through $v$, so that the operator consists of the usual first order partial derivatives in $y$ and $v$, and the fractional derivatives are in terms of $y$ only. Then, using (\ref{eqn:RieszFellerDeriv}), we obtain the corresponding partial differential equation for the CF,
\begin{equation}
\begin{cases}
&\displaystyle \frac{\dd \psi}{\dd t} = im\mathcal{F}\left[\hat{f}(\bar{x},y)p\right] + \frac{il}{\epsilon}\mathcal{F}\left[ g(\bar{x},y)p\right] - \frac{b^{\alpha}|l|^{\alpha}}{\epsilon^{(1-\gamma)\alpha}}\Xi(l;\alpha,\beta)\psi, \\ &\psi(l,m,0) = \exp(ily_{0}).
\end{cases}  \label{FT_FPE_L_canon_reduced}
\end{equation}
Substituting $g$ and $\hat{f}$ as in (\ref{drifts}) and (\ref{v_dynamics}) into (\ref{FT_FPE_L_canon_reduced}) yields
\begin{align}
\begin{cases}
&\displaystyle \frac{\dd \psi}{\dd t} = \left( \frac{f_{2}}{\epsilon^{\gamma}}m + \frac{g_{2}}{\epsilon}l\right)\frac{\dd \psi}{\dd l} + \left[i\frac{g_{1}}{\epsilon^{1-\gamma}}l + i\frac{g_{1}f_{2}}{g_{2}}m - \frac{b^{\alpha}|l|^{\alpha}}{\epsilon^{(1-\gamma)\alpha}}\Xi(l;\alpha,\beta)\right]\psi,\\ &\psi(l,m,0) = \exp(ily_{0}) \end{cases} \label{FT_FPE_L_gener}
\end{align}
where we have dropped the argument $\bar{x}$ from $f_{1},\,f_{2},\,g_{1},\,g_{2}$. We derive the solution for $\psi$ using the method of characteristics. The characteristic curves for the variables $t,l$ and $m$ are given in terms of the characteristic variables $r,l_{0}$, and $m_{0}$,
\begin{equation}
t(r) = r, \quad l(r) = l_{0}\exp\left(-\frac{g_{2}}{\epsilon}r \right) - \epsilon^{1-\gamma}\frac{f_{2}}{g_{2}}m_{0}\left( 1 - \exp\left(-\frac{g_{2}}{\epsilon}r \right) \right), \quad m(r) = m_{0},
\end{equation}
and $\psi$ satsifies
\begin{equation}
\frac{d\psi}{dr} = \left(i\frac{f_{2}g_{1}}{g_{2}}m + i\epsilon^{\gamma-1}g_{1}l -\frac{b^{\alpha}}{\epsilon^{(1-\gamma)\alpha}}|l(r)|^{\alpha}\Xi(l(r);\alpha,\beta)\right)\psi, \quad \psi(0) = \exp(il_{0}y_{0}) \label{method_of_chars_psi}
\end{equation}
where $\Xi$ is given in (\ref{Xi}). The solution of (\ref{method_of_chars_psi}) is given by
\begin{align}
\psi(l,m,t) = \exp&\left[i\Lambda(t)y_{0} - i\left(\epsilon^{\gamma}\frac{g_{1}}{g_{2}}l + \epsilon\frac{g_{1}f_{2}}{g_{2}^{2}}m\right)\left(1 - \exp(g_{2}t/\epsilon)\right) \right. \label{psi_integral} \\ & \qquad \qquad \left. -\frac{b^{\alpha}}{\epsilon^{(1-\gamma)\alpha}}\int_{0}^{t}|\Lambda(r)|^{\alpha}\Xi(\Lambda(r);\alpha,\beta)\,dr \right], \nonumber 
\end{align}
where
\begin{equation}
\Lambda(s) = l\exp\left(\frac{g_{2}}{\epsilon}s\right) - \epsilon^{1-\gamma}m\frac{f_{2}}{g_{2}}\left(1 - \exp\left(\frac{g_{2}}{\epsilon}s\right)\right). \label{Lambda}
\end{equation}
We evaluate the integral term of (\ref{psi_integral}) asymptotically for small $\epsilon$ as shown in Appendix \ref{app:calculations} giving the result,
\begin{align}
&\int_{0}^{t}|\Lambda(r)|^{\alpha}\Xi(\Lambda(r);\alpha,\beta)\,dr \label{integral_result} \\  & = \begin{cases}\displaystyle \left(\frac{\epsilon|l|^{\alpha}}{-\alpha g_{2}}\right)\Xi(l;\alpha,\beta)(1 - \exp({\alpha g_{2} t / \epsilon})) + O(\epsilon^{2 - \gamma}) &\mbox{for $t = O(\epsilon),$} \\ \displaystyle  \epsilon|f_{2}g_{2}^{-1}|^{\alpha}|m|^{\alpha}\Xi(m;\alpha,\beta^{*})t + \frac{\epsilon|l|^{\alpha}}{-\alpha g_{2}}\Xi(l;\alpha,\beta) + O(\epsilon^{2-\gamma}) &\mbox{for $t = O(1)$.} \end{cases} \nonumber
\end{align}
where 
\begin{equation}
\beta^{*} = \beta\sgn{f_{2}}. \label{beta_star}
\end{equation}
Substituting (\ref{integral_result}) into (\ref{psi_integral}) gives the following approximation for $\psi$,
\begin{equation}
\psi \sim \begin{cases}\exp\left[i\Lambda(t)y_{0} - i\left(\epsilon^{\gamma}\frac{g_{1}}{g_{2}}l + \epsilon\frac{g_{1}f_{2}}{g_{2}^{2}}m\right)\left(1 - \exp(g_{2}t/\epsilon)\right) \right. \\ \qquad \qquad \left. -\frac{b^{\alpha}|l|^{\alpha}}{\epsilon^{(1-\gamma)\alpha - 1}(-\alpha g_{2})}\Xi(l;\alpha,\beta)(1 - \exp({\alpha g_{2} t / \epsilon})) \right] \quad &\mbox{for $t = O(\epsilon)$,} \\ \exp\left[ -il\epsilon^{\gamma}\frac{g_{1}}{g_{2}} - im\left(\epsilon\frac{g_{1}f_{2}}{g_{2}^{2}} + \epsilon^{1-\gamma}\frac{f_{2}}{g_{2}}y_{0}\right) -\frac{b^{\alpha}|l|^{\alpha}}{\epsilon^{(1-\gamma)\alpha-1}(-\alpha g_{2})}\Xi(l;\alpha,\beta) \right. \\ \qquad \qquad \left. -\frac{b^{\alpha}}{\epsilon^{(1-\gamma)\alpha-1}}|f_{2}/g_{2}|^{\alpha}|m|^{\alpha}\Xi(m;\alpha,\beta^{*})t   \right] \quad &\mbox{for $t = O(1).$} \end{cases} \label{psi_approx}
\end{equation}
Since $\psi$ has no terms that involve products of $l$ and $m$, we can factor the approximate CF into components $\psi \approx \psi_{y}{\psi_{v}}$ that correspond to CFs of $y_{t}$ and $v_{t}$. The CF for $y_{t}$ given in (\ref{canon2}) and $g$ given in (\ref{drifts}), treating $\bar{x}_{t}$ as a constant is,
\begin{align}
\psi_{y}(l,t) = \exp&\left[ ily_{0}\exp\left( {g_{2}t/\epsilon}\right) -il\epsilon^{\gamma}\frac{g_{1}}{g_{2}}(1 - \exp\left( {g_{2}t}/{\epsilon}\right)) \right.\label{char_func_y} \\  &\qquad \qquad \left.-\frac{b^{\alpha}|l|^{\alpha}}{\epsilon^{(1-\gamma)\alpha-1}(-\alpha g_{2})}\Xi(l;\alpha,\beta)(1 - \exp({\alpha g_{2} t / \epsilon}))\right]. \nonumber
\end{align}
The asymptotic approximation for $t = O(1)$ is
\begin{align}
\psi_{y}(l,t) &\sim \exp\left[ -il\epsilon^{\gamma}\frac{g_{1}}{g_{2}}  -\frac{b^{\alpha}|l|^{\alpha}}{\epsilon^{(1-\gamma)\alpha-1}(-\alpha g_{2})}\Xi(l;\alpha,\beta)\right]. \label{char_func_y_asymp}
\end{align}
We then factor $\psi_{y}$ as defined by (\ref{char_func_y}) for $t = O(\epsilon)$ and (\ref{char_func_y_asymp}) for $t = O(1)$ out of the approximate expression for $\psi$ given in (\ref{psi_approx}) to get the asymptotic behaviour of $\psi_{v}$,
\begin{align}
\psi_{v}(m,t) &\sim \exp\left[-im\left(\epsilon^{1-\gamma}\frac{f_{2}}{g_{2}}y_{0} + \epsilon\frac{g_{1}f_{2}}{g_{2}^{2}}\right)(1 - \exp(g_{2}t/\epsilon))\right] &\mbox{for $t = O(\epsilon)$,} \label{Capital_Psi_small_t} \\ \psi_{v}(m,t) &\sim \exp\left[ -im\left(\epsilon^{1-\gamma}\frac{f_{2}}{g_{2}}y_{0} + \epsilon\frac{g_{1}f_{2}}{g_{2}^{2}}\right) \right. \label{Capital_Psi} \\ & \qquad \qquad \left. - \frac{b^{\alpha}}{\epsilon^{(1-\gamma)\alpha - 1}}|f_{2}g_{2}^{-1}|^{\alpha}|m|^{\alpha}\Xi(m;\alpha,\beta^{*})t   \right] &\mbox{for $t = O(1).$} \nonumber
\end{align} 
This factorization of the CF allows us to analyze the asymptotic approximation of the CF, $\psi_{v}$ for the process $v_{t}$. We compare (\ref{Capital_Psi}) to the CF of an $\alpha$-stable process, $z_{t},\, t \ge 0$, with the SDE and associated CF given by
\begin{equation}
dz_{t} = \Sigma\,dL_{t}^{(\alpha,\beta)}, \quad \psi_{z}(k) = \exp\left[ikz_{0} - t\Sigma^{\alpha}|k|^{\alpha}\Xi(k;\alpha,\beta)\right], \label{char_WLp}
\end{equation} 
where $\Sigma$ is a constant and $k$ is the Fourier variable. The functional structure with respect to the Fourier variables $m$ in (\ref{Capital_Psi}) and $k$ in (\ref{char_WLp}) is identical, with the terms $(\alpha,\beta,\Sigma)$ in (\ref{char_WLp}) corresponding to $(\alpha,\beta^{*}, \epsilon^{(\gamma - 1 + 1/\alpha)}b|f_{2}/g_{2}|)$ in (\ref{Capital_Psi}). The coefficient of $im$ in (\ref{Capital_Psi_small_t}), (\ref{Capital_Psi}) defines the mean of $v_{t}$, which is non-zero in general if $y_{0} \neq \expectedcond{y}{y|x}$ given in (\ref{fbar_specific}). The mean of $v_{t}$ is constant on the $t = O(1)$ time scale as is evident by inspection of (\ref{Capital_Psi}) and is thus inconsequential to our analysis since we are interested in the increments of $v_{t}$ on this time scale rather than $v_{t}$ itself.

Since convergence in CF implies convergence in distribution \cite{Feller1966II}, the comparison of (\ref{Capital_Psi}) and (\ref{char_WLp}) indicates that on the $O(1)$ time scale, the increments $dv_{t}$ are approximately distributed according to an $\alpha$-stable law, $dv_{t} \sim \epsilon^{\rho}b|f_{2}/g_{2}|\,dL_{t}^{(\alpha,\beta^{*})}$ where $\rho = \gamma - 1 + 1/\alpha$. We conclude that the (L) approximation for the system (\ref{canon1} -- \ref{drifts}) is given by $\bar{x}_{t} + \xi_{t}$ where $\bar{x}_{t}$ is given by (\ref{Aapprox}) and
\begin{equation}
d\xi_{t} = \bar{f}'(\bar{x}_{t})\xi_{t}\,dt + \epsilon^{\rho}b\left|\frac{f_{2}(\bar{x}_{t})}{g_{2}(\bar{x}_{t})}\right|\,dL_{t}^{(\alpha,\beta^{*})}, \quad \xi_{0} = 0 \label{Lapprox}
\end{equation}
where $\bar{f}$ is given in (\ref{fbar_specific}). We note that if $\rho < 0$ (i.e. $\gamma < 1 - 1/\alpha$), then the scale parameter of the noise diverges and the approximation (\ref{Lapprox}) is undefined in the $\epsilon \rightarrow 0$ limit. If $\rho > 0$, then the scaling of the noise approaches zero in this limit and $\xi_{t} \rightarrow 0$ as $t \rightarrow \infty$, provided $\bar{f}'(\bar{x}_{t}) < 0$ for all $\bar{x}_{t}$. In this case, the (L) approximation is asymptotically equivalent to the (A) approximation (\ref{Aapprox}) as $\epsilon \rightarrow 0$. For non-zero, but sufficiently small $\epsilon$, the (L) approximation is an asymptotic approximation in the weak sense for the dynamics of $x_{t}$ for $\rho > 0$. This is also true for $\rho < 0$, $\epsilon$ not too small. However in this case, the dynamics of $x_{t}$ become dominated by the noise and it becomes questionable that $x_{t}$ can be considered a slow process relative to $y_{t}$. For $\rho = 0$, the intensity of the noise in (\ref{Lapprox}) does not depend on $\epsilon$ and the stochastic corrections persist as the timescale separation is made arbitrarily small. 

The (L) approximation is useful not only when the system being approximated is itself linear, but also when the system is near a deterministic attractor to which $\bar{x}_{t}$ has converged. In this case, the (L) approximation is reasonable on finite time intervals that are shorter than the escape time of the local attractor. We compare simulations of the full system and the (L) approximation for both linear and nonlinear systems in \S \ref{sec:comp_results}.

\subsection{(N+) approximation}
\label{subsec:N+approx}

The (L) approximation is a good approximation for linear systems or for certain local approximations. However, for systems with nonlinear dynamics, the (L) approximation may be inadequate due to its inability to model systems with multimodal stationary distributions or behaviour different from an OULp. In these cases, the (N+) approximation is more appropriate as it is a better approximation for such nonlinear systems.

Beginning with the original system (\ref{canon1}, \ref{canon2}) with $f,g$ given by (\ref{drifts}), we express the dynamics of $x_{t}$ in terms of $\bar{f}$ and $\hat{f}$ given in (\ref{fbar}) and (\ref{v_dynamics}),
\begin{equation}
dx_{t} = \bar{f}(x_{t})\,dt + \hat{f}(x_{t},y_{t})\,dt.
\end{equation}
As in the previous section, to complete the (N+) approximation, we find a stochastic approximation for $\hat{f}(x_{t},y_{t})\,dt$ that has its distributional properties on the $t = O(1)$ time scale. We make a change of coordinates $\mu_{t} = \mathcal{T}(x_{t})$ where $\mathcal{T}$ is a differentiable function satisfying
\begin{equation}
\mathcal{T}'(x)f_{2}(x) = g_{2}(x). \label{eqn:Ttransformation}
\end{equation}
The transformation $\mu_{t} = \mathcal{T}(x_{t})$ defined in (\ref{eqn:Ttransformation}) is invertible since $g_{2}(x) < 0$ and $f_{2}(x) \ne 0$ has the same sign for all $x$. As we show below, this coordinate transformation, $x \rightarrow \mu$ allows us to approximate the increments $\hat{f}(x_{t},y_{t})\,dt$ in the transformed coordinates as an additive $\alpha$-stable noise process when we consider the CF in Fourier space. Then there is no ambiguity in the interpretation of these noise terms, particularly when we invert the transformation $\mathcal{T}$ to obtain the approximation in the original variables. Our use of such a coordinate transformation (\ref{eqn:Ttransformation}) is similar to that in \cite{Chechkin2014} where a fast-slow system is transformed so that the influence of the fast Ornstein-Uhlenbeck process drives the slow variable additively. In \cite{Chechkin2014}, this formulation allows for the unambiguous derivation of a reduced stochastic model, with a Stratonovich interpretation. 

The dynamics of $\mu_{t}$ can be determined using the change of variables formula \cite{Cont2004},
\begin{align} 
d\mu_{t} &= \mathcal{T}'(x_{t-})\,dx_{t} + \frac{1}{2}\mathcal{T}''(x_{t-})\,d[x,x]_{t}^{c} + \mathcal{T}(x_{t}) - \mathcal{T}(x_{t-}) - \mathcal{T}'(x_{t})\Delta x_{t},
\end{align}
where $x_{t-} = \lim_{s \rightarrow t^{-}}x_{s}$ and $\Delta x_{t}$ denotes the jump component of $x_{t}$. Since the continuous part of the quadratic variation of $x$, $d[x,x]_{t}^{c} \propto (dt)^{2} \rightarrow 0$ and there are no jump components because $x_{t}$ is continuous for all $t$, we get
\begin{align}
d\mu_{t} &= \mathcal{T}'(x_{t})\,dx_{t} \\ 
&= \mathcal{T}'(\mathcal{T}^{-1}(\mu_{t}))\left[ \bar{f}(\mathcal{T}^{-1}(\mu_{t})) + \hat{f}(\mathcal{T}^{-1}(\mu_{t}),y_{t}) \right]\,dt \\ 
&= \mathcal{T}'(\mathcal{T}^{-1}(\mu_{t}))\bar{f}(\mathcal{T}^{-1}(\mu_{t}))\,dt + \left[ g_{1}(\mathcal{T}^{-1}(\mu_{t})) + \epsilon^{-\gamma}g_{2}(\mathcal{T}^{-1}(\mu_{t}))y_{t}\right]\,dt. \label{eqn:mu_dynamics}
\end{align}
We define $dV_{t}$ as the additional fluctuations about the first forcing term in (\ref{eqn:mu_dynamics}), analogous to $dv_{t}$ in the (L) approximation (\ref{v_dynamics}),
\begin{equation}
dV_{t} = \left[ g_{1}(\mathcal{T}^{-1}(\mu_{t})) + \epsilon^{-\gamma}g_{2}(\mathcal{T}^{-1}(\mu_{t}))y_{t}\right]\,dt, \quad V_{0} = 0. \label{eqn:v_dynamics_N+1}
\end{equation}
Note that the dynamics in (\ref{eqn:v_dynamics_N+1}) are the same (up to a multiplicative constant) as the drift terms in the equation for $y$. This correspondence between the $V$ and $y$ equations leads to a straightforward factorization of the joint CF in Fourier space, as we see below. 

Invoking the separation of time scales for small $\epsilon$, we use a multiple scale approximation for the FKE for the $(y,V)$ system treating $\mu_{t}$ as effectively constant relative to the fast time scales of $y_{t}$ and $V_{t}$,
\begin{equation}
\begin{cases}
&\displaystyle \frac{\dd q}{\dd t} = -\left(\frac{\dd}{\dd V} + \epsilon^{\gamma-1}\frac{\dd}{\dd y} \right)\left(\left[ g_{1} + \epsilon^{-\gamma}g_{2}y\right]q \right) + \frac{b^{\alpha}}{\epsilon^{(1 - \gamma)\alpha}}\mathcal{D}_{y}^{(\alpha,\beta)}q, \\ &q(y,V,0) = \delta(y - y_{0})\delta(V) \end{cases} \label{eqn:FKE_vw}
\end{equation}
where $q = q(y,V,t)$ is the joint PDF for $y,V$ and $g_{1}$ and $g_{2}$ are shorthand for $g_{1}(\mathcal{T}^{-1}(\mu_{t}))$ and $g_{2}(\mathcal{T}^{-1}(\mu_{t}))$, treated as constants. We take the FT of (\ref{eqn:FKE_vw}) and define the CF of $(y,V)$ to be $\phi = \phi(l,m,t)$ where
\begin{equation}
\phi(l,m,t) = \iint_{\mathbb{R}^{2}}\exp(ily + imV)q(y,V,t)\,dV dy. \label{char_func2}
\end{equation}
The resulting linear partial differential equation for $\phi$ is similar to (\ref{FT_FPE_L_gener}),
\begin{align}
\begin{cases}& \displaystyle \frac{\dd \phi}{\dd t} = \left( \frac{g_{2}}{\epsilon^{\gamma}}m + \frac{g_{2}}{\epsilon}l\right)\frac{\dd \phi}{\dd l} + \left[i\left(g_{1}m + \frac{g_{1}l}{\epsilon^{1-\gamma}}\right) - \frac{b^{\alpha}|l|^{\alpha}}{\epsilon^{(1-\gamma)\alpha}}\Xi(l;\alpha,\beta)\right]\phi, \\  &\phi(l,m,0) = \exp(ily_{0}) \end{cases} \label{FT_FPE_N_gener}
\end{align}
with $\Xi(l;\alpha,\beta)$ given in (\ref{Xi}). Analogous to (\ref{FT_FPE_L_gener}), (\ref{FT_FPE_N_gener}) can be solved using the method of characteristics giving us
\begin{align}
\phi(l,m,t) = \exp&\left[ i\Gamma(t)y_{0} - i\frac{g_{1}}{g_{2}}(\epsilon^{\gamma}l + \epsilon m)(1 - \exp(g_{2}t/\epsilon)) \right. \label{phi_integral} \\  & \qquad\qquad \left.-\frac{b^{\alpha}}{\epsilon^{(1-\gamma)\alpha}}\int_{0}^{t}|\Gamma(r)|^{\alpha}\Xi(\Gamma(r);\alpha,\beta)\,dr \right], \nonumber
\end{align}
where
\begin{equation}
\Gamma(s) = l\exp\left(\frac{g_{2}}{\epsilon}s\right) - \epsilon^{1-\gamma}m\left(1 - \exp\left(\frac{g_{2}}{\epsilon}s\right)\right). \label{Gamma}
\end{equation}
Using the same asymptotic results derived in Appendix \ref{app:calculations}, we can approximate $\phi$ on the $t = O(\epsilon)$ and $O(1)$ time scales. Thus, we obtain an approximate form for the CF of the $(y,V)$ system for $0 < \epsilon \ll 1$,
\begin{equation}
\phi \sim \begin{cases}\exp\left[i\Gamma(t)y_{0} - i(\epsilon^{\gamma}l + \epsilon m)\frac{g_{1}}{g_{2}}(1 - \exp(g_{2}t/\epsilon)) \right. \\ \qquad \left. -\frac{b^{\alpha}|l|^{\alpha}\Xi(l;\alpha,\beta)}{\epsilon^{(1-\gamma)\alpha-1}(-\alpha g_{2})}(1 - \exp({\alpha g_{2} t / \epsilon})) \right] &\mbox{for $t = O(\epsilon)$} \\ 
\exp\left[ - i\epsilon^{\gamma}\frac{g_{1}}{g_{2}}l - i\left(\epsilon^{1-\gamma}y_{0} + \epsilon\frac{g_{1}}{g_{2}} \right)m \right. \\ \left. \qquad -\frac{b^{\alpha}|l|^{\alpha}\Xi(l;\alpha,\beta)}{\epsilon^{(1-\gamma)\alpha-1}(-\alpha g_{2})} - \frac{b^{\alpha}}{\epsilon^{(1-\gamma)\alpha - 1}}|m|^{\alpha}\Xi(m;\alpha,-\beta)t\right] &\mbox{for $t = O(1).$} \end{cases} \label{phi_approx}
\end{equation}
Since $\phi$ has no products of $l$ and $m$, we decompose the $t = O(1)$ approximation in (\ref{phi_approx}) into two factors: $\phi \sim \phi_{y}\phi_{V}$, where 
\begin{equation}
\phi_{V}(m,t) = \exp\left(- i\left(\epsilon^{1-\gamma}y_{0} + \epsilon\frac{g_{1}}{g_{2}} \right)m - \epsilon^{\alpha\rho}b^{\alpha}|m|^{\alpha}\Xi(m;\alpha,-\beta)t\right)
\end{equation}
and $\phi_{y} = \psi_{y}$ is the CF of $y$ given in (\ref{char_func_y_asymp}). As in the case of the (L) approximation, $\phi_{V}$ is identical in form to the CF of an $\alpha$-stable process. By the same reasoning leading from (\ref{Capital_Psi}) to (\ref{Lapprox}) in the previous section, we conclude that the dynamics of the perturbation $V_{t}$ can be approximated on the $t = O(1)$ time scale by an $\alpha$-stable noise process,
\begin{equation}
dV_{t} =  g_{1}(\mathcal{T}^{-1}(\mu_{t}))\,dt + \epsilon^{-\gamma}g_{2}(\mathcal{T}^{-1}(\mu_{t}))y_{t}\,dt \approx \epsilon^{\rho}b\,dL_{t}^{(\alpha,-\beta)}. \label{eqn:v_dynamics_N+2}
\end{equation}
Note that in our transformed system $\mu_{t} = \mathcal{T}(x_{t})$, the coefficient of the L\'{e}vy increment is state-independent and so the limiting noise is additive. Using (\ref{eqn:v_dynamics_N+1}) and (\ref{eqn:v_dynamics_N+2}) in (\ref{eqn:mu_dynamics}), we obtain an approximate equation for $\mu_{t}$,
\begin{equation}
d\mu_{t} \approx \mathcal{T}'(\mathcal{T}^{-1}(\mu_{t}))\bar{f}(\mathcal{T}^{-1}(\mu_{t}))\,dt - \epsilon^{\rho}b\,dL_{t}^{(\alpha,\beta)}, \quad \mu_{0} = \mathcal{T}(x_{0}).
\end{equation}
Here we have chosen to express the stochastic increments $\epsilon^{\rho}b\,dL_{t}^{(\alpha,-\beta)}$ as $-\epsilon^{\rho}b\,dL_{t}^{(\alpha,\beta)}$, which is equivalently distributed. To complete the approximation for $x_{t}$ on the slow time scale, we invert our initial transformation $x_{t} \rightarrow \mu_{t} = \mathcal{T}(x_{t})$ by taking ${\mu_{t}} \rightarrow X_{t} = \mathcal{T}^{-1}({\mu_{t}})$, similarly to \cite{Chechkin2014} in the example of Gaussian noise. The fact that the stochastic forcing on $\mu_{t}$ is additive allows for an unambiguous interpretation of the noise terms prior to performing the inverse transformation, and inverting the transformation then provides the appropriate interpretation of the noise terms for the approximation of $x_{t}$. Since $\mathcal{T}^{-1}$ is a nonlinear transformation, the change of variables formula for jump processes \cite{Cont2004} yields non-trivial contributions when applied to $X_{t}$ as follows,
\begin{align}
d{X}_{t} &= (\mathcal{T}^{-1})'(\mu_{t-})\,d\mu_{t} + \mathcal{T}^{-1}(\mu_{t})  - \mathcal{T}^{-1}(\mu_{t-}) - (\mathcal{T}^{-1})'(\mu_{t-})\,\Delta\mu_{t}, \label{eqn:xhatdynamics_intermed}
\end{align}
where $X_{0} = x_{0}$, $\mu_{t-} = \lim_{s \rightarrow t^{-}}\mu_{s}$, and $\Delta\mu_{t} = \mu_{t} - \mu_{t-} = -\epsilon^{\rho} b\,dL_{t}^{(\alpha,\beta)}$ is the jump component of $\mu_{t}$ at time $t$. Using the expression for $d\mu_{t}$ (\ref{eqn:mu_dynamics}), (\ref{eqn:xhatdynamics_intermed}) can be written as 
\begin{align}
d{X}_{t} &= \frac{1}{\mathcal{T}'(\mathcal{T}^{-1}(\mu_{t-}))}\mathcal{T}'(\mathcal{T}^{-1}(\mu_{t}))\bar{f}(\mathcal{T}^{-1}(\mu_{t}))\,dt + \left\{\mathcal{T}^{-1}(\mu_{t})  - \mathcal{T}^{-1}(\mu_{t-})\right\} \\
&= \left(\frac{\mathcal{T}'({X}_{t})}{\mathcal{T}'({X}_{t-})}\right)\bar{f}({X}_{t})\,dt + \left\{\mathcal{T}^{-1}(\mu_{t})  - \mathcal{T}^{-1}(\mu_{t-})\right\}. \label{eqn:xhatdynamics_N+}
\end{align}
where $X_{t-} = \lim_{s \rightarrow t^{-}}X_{s}$. With probability one, the ratio $\mathcal{T}'({X}_{t}) / \mathcal{T}'({X}_{t-})$ is equal to one, as ${X}_{t}$ is a c\`{a}dl\`{a}g process and hence the set of times where a jump occurs is a null-set with respect to the Lebesgue measure \cite{Cont2004}. Hence, we replace $\mathcal{T}'({X}_{t}) / \mathcal{T}'({X}_{t-})$ with $1$ in (\ref{eqn:xhatdynamics_N+}). The difference $\mathcal{T}^{-1}(\mu_{t})  - \mathcal{T}^{-1}(\mu_{t-})$ can not be neglected since the time integral of these terms results in a sum that does not scale with the size of the infinitesimal time step $dt$. To see this, we show that this difference is equal to the Marcus increment $[\mathcal{T}'({X}_{t})]^{-1} \diamond dL_{t}^{(\alpha,\beta)}$, by writing $\mathcal{T}^{-1}(\mu_{t})  - \mathcal{T}^{-1}(\mu_{t-})$ in terms of the processes, $\Theta(r; \Delta L_{t}, \mu_{t-}) = \mu_{t-} - \epsilon^{\rho}b\Delta L_{t}r$ and $\theta(r;\Delta L_{t}, \mathcal{T}^{-1}(\mu_{t-})) = \mathcal{T}^{-1}(\Theta(r;\Delta L_{t}, \mu_{t-}))$, where
\begin{equation}
\frac{d\theta}{dr} = \frac{d\theta}{d\Theta}\frac{d\Theta}{dr} = -\frac{\epsilon^{\rho}b\Delta L_{t}}{\mathcal{T}'(\theta(r))}, \quad \theta(0; \Delta L_{t}, \mathcal{T}^{-1}(\mu_{t-})) = \mathcal{T}^{-1}(\mu_{t-}) = {X}_{t-} \label{eqn:Marcus_N+2}
\end{equation}
and $\Delta L_{t} = dL_{t}^{(\alpha,\beta)}$ is the size of the jump in the $\alpha$-stable noise process $L_{t}^{(\alpha,\beta)}$ at time $t$. The variable $\theta(r; \Delta L_{t},\mu_{t-})$ is referred to as the Marcus map and it describes the integrated effect of the jump $\Delta L_{t}$ over the infinitesimal time on which the jump occurs which begins at $r = 0$ and ends at $r = 1$ (additional details in Appendix \ref{app:simulation_Marcus}). The difference $\mathcal{T}^{-1}(\mu_{t})  - \mathcal{T}^{-1}(\mu_{t-})$ can then be written in terms of $\theta(r; \Delta L_{t},\mu_{t-})$ as follows,
\begin{align}
&\mu_{t} = \mu_{t-} - \epsilon^{\rho}b\Delta L_{t} = \Theta(1; \Delta L_{t}, \mu_{t-}) = \mathcal{T}(\theta(1; \Delta L_{t}, \mathcal{T}^{-1}(\mu_{t-}))) \\
\Rightarrow \quad &\mathcal{T}^{-1}(\mu_{t})  - \mathcal{T}^{-1}(\mu_{t-}) = \mathcal{T}^{-1}(\Theta(1; \Delta L_{t}, \mu_{t-}))  - X_{t-} = \theta(1; \Delta L_{t}, X_{t-}) - X_{t-} \label{eqn:intermediate_marcus_example}
\end{align}
The difference $\theta(1; \Delta L_{t}, X_{t-}) - X_{t-}$ where $\theta$ satisfies (\ref{eqn:Marcus_N+2}) is the definition of the Marcus increment $-\epsilon^{\rho}b[\mathcal{T}'({X}_{t})]^{-1}\diamond dL_{t}^{(\alpha,\beta)}$. Then using (\ref{eqn:intermediate_marcus_example}) and $\mathcal{T}'$ from (\ref{eqn:Ttransformation}) in (\ref{eqn:xhatdynamics_N+}), yields the weak approximation $X_{t}$ for the dynamics of $x_{t}$,
\begin{equation}
d{X}_{t} = \bar{f}({X}_{t})\,dt + \epsilon^{\rho}b\left(\frac{f_{2}({X}_{t})}{-g_{2}({X}_{t})}\right)\diamond dL_{t}^{(\alpha,\beta)}, \quad X_{0} = x_{0}. \label{N+approx}
\end{equation}
We note that (\ref{N+approx}), which we refer to as the (N+) approximation is equivalent to the (L) approximation (\ref{Lapprox}) when (\ref{canon1}, \ref{canon2}) are a linear system of equations (i.e. $f_{1},f_{2},g_{1},g_{2}$ are constant). Also note that the stochastic increment $-\epsilon^{\rho} b [f_{2}({X}_{t})/g_{2}({X}_{t})]\diamond dL_{t}^{(\alpha,\beta)}$ can be written equivalently as $\epsilon^{\rho}b\left|f_{2}({X}_{t})/(g_{2}({X}_{t})) \right|\diamond dL_{t}^{(\alpha,\beta^{*})}$ where $\beta^{*}$ is given in (\ref{beta_star}) with $f_{2} = f_{2}({X}_{t})$ since any negative sign in the scale parameter can be absorbed into the skewness parameter. As in \S \ref{subsec:Lapprox}, if $\rho > 0$, then the (N+) approximation reduces to the (A) approximation (\ref{Aapprox}) in the limit $\epsilon \rightarrow 0$, while the scale parameter diverges if $\rho < 0$ in the same limit. In the case of small, non-zero $\epsilon$, (\ref{N+approx}) provides a good approximation to $x_{t}$ when $\rho > 0$, but if $\rho < 0$ then the approximation may not be appropriate for the same reasons given in \S \ref{subsec:Lapprox}. If $\rho = 0$, the stochastic dynamics persist in the strict limit $\epsilon \rightarrow 0$.

When $\alpha = 2$, the Marcus interpretation of (\ref{N+approx}) reduces to the Stratonovich interpretation (see Appendix \ref{app:simulation_Marcus}) and thus our approximation is consistent with the previously derived (N+) approximation in the case of Gaussian white noise stochastic forcing \cite{Arnold2003}.

In the next section we show that the (N+) approximation is generally superior to the (L) approximation for nonlinear systems as would be expected. However, we note that the (N+) approximation requires an additional constraint (i.e. $\sgn{f_{2}(x)}$ is constant and not equal to zero, for all $x$) and its numerical implementation is more complicated, often requiring simulation of the Marcus differential terms (See Appendix \ref{app:simulation_Marcus}) for which we do not have an efficient numerical method when the solution for $\theta$ from (\ref{eqn:Marcus_N+2}) cannot be expressed in closed form and  $\beta \ne 0$. Depending on the problem under consideration, for practical reasons the (N+) approximation may not be a preferable option over the (L) approximation. 

Finally, we note the interesting fact that the drift and diffusion coefficients of the dynamics in the approximations are identical to those in the Gaussian case discussed in \cite{Mona11} up to a power of $\epsilon$, regardless of the value of $\alpha$.

\section{Computational Results}
\label{sec:comp_results}
We apply the (L) and (N+) stochastic averaging approximations to four fast-slow dynamical systems with $\alpha$-stable noise forcing and numerically simulate the systems to compare the stationary probability density functions and autocodifference functions of the approximations to those of the corresponding full systems. One of the systems we simulate is linear while the other three include nonlinear functions $f_{1}, f_{2}, g_{1}$ and $g_{2}$. For this section, we set $\rho = 0$ (i.e. $\gamma = 1 - 1/\alpha$). Our numerical scheme is outlined in Appendix \ref{app:simulating_SDEs}, and we have run additional simulations with even smaller simulation time steps to confirm that the simulations are well-resolved. All density plots in this paper were derived from simulations having the same sampling time step $Dt = 10^{-2}$ and number of data points, $N = 10^{8}$. The rationale for our choice of time series length is given in Appendix \ref{app:stats}.

\subsection{Linear System}
\label{subsec:linsys}
We apply the (L) approximation from \S \ref{subsec:Lapprox} to a linear system inspired by an example in \cite{Mona11},
\begin{align}
dx_{t} &= \left(-x_{t} + \epsilon^{-\gamma}ay_{t} \right)\,dt, \quad x_{0} \in \mathbb{R} \label{linsys1}\\
\epsilon\,dy_{t} &= \left(\epsilon^{\gamma}cx_{t} - y_{t} \right)\,dt + \epsilon^{\gamma}b\,dL_{t}^{(\alpha,\beta)}, \quad y_{0} = 0 \label{linsys2}
\end{align}
where $0 < \epsilon \ll 1$ and $(1 - ac) > 0$. Using (\ref{Aapprox}), (\ref{Lapprox}), we write the (L) approximation for $x_{t}$ in (\ref{linsys1}, \ref{linsys2}). The (A) approximation, $\bar{x}$, as in (\ref{Aapprox}) is
\begin{equation}
d\bar{x}_{t} = -(1 - ac)\bar{x}_{t}\,dt, \quad \bar{x}_{0} = x_{0}, \quad \Rightarrow \quad \bar{x}_{t} = x_{0}\exp[-(1 - ac)t]. \label{Aapprox_lin}
\end{equation}
Using (\ref{Lapprox}) and (\ref{Aapprox_lin}), the (L) approximation for $x_{t}$ is $\bar{x}_{t} + \xi_{t}$ where $\bar{x}_{t}$ satisfies (\ref{Aapprox_lin}) and $\xi_{t}$ satisfies
\begin{equation}
d\xi_{t} = -(1 - ac)\xi_{t}\,dt + |a|b\,dL_{t}^{(\alpha,\beta_{a})}, \quad \xi_{0} = 0, \label{Lapprox_linear}
\end{equation}
where $\beta_{a} = \sgn{a}\beta$. In the long-time limit $t \rightarrow \infty$ (i.e. at stationarity) or in the case that $x_{0} = 0$, the (N+) approximation of $x_{t}$ is identical to the (L) approximation since (\ref{linsys1}, \ref{linsys2}) is linear. We set $x_{0} = 0$ for all the numerical simulations of (\ref{linsys1}, \ref{linsys2}) and (\ref{Lapprox_linear}), and hence $\bar{x}_{t} = 0$.

\begin{figure}[h]
\centering
\includegraphics[width=0.323\textwidth]{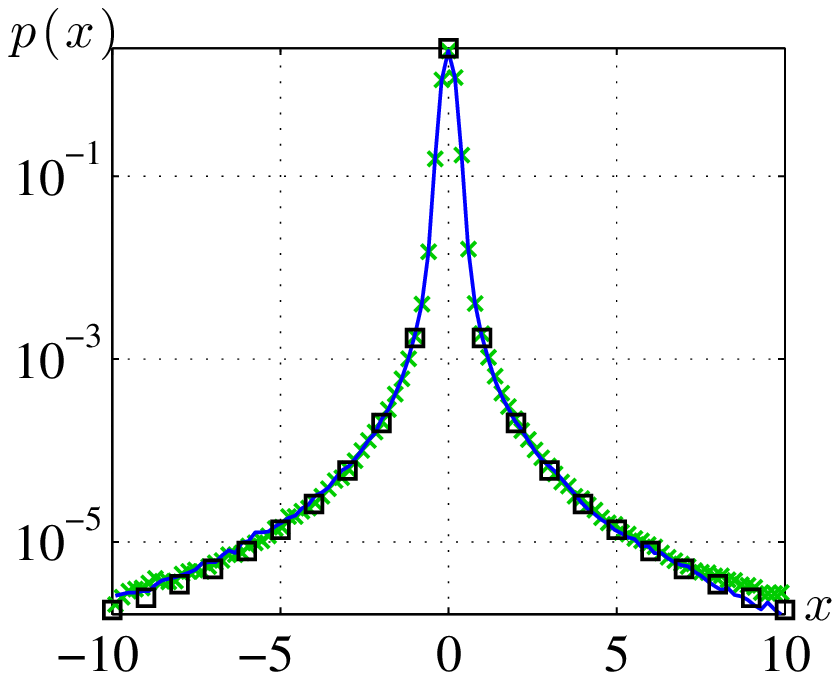}
\includegraphics[width=0.323\textwidth]{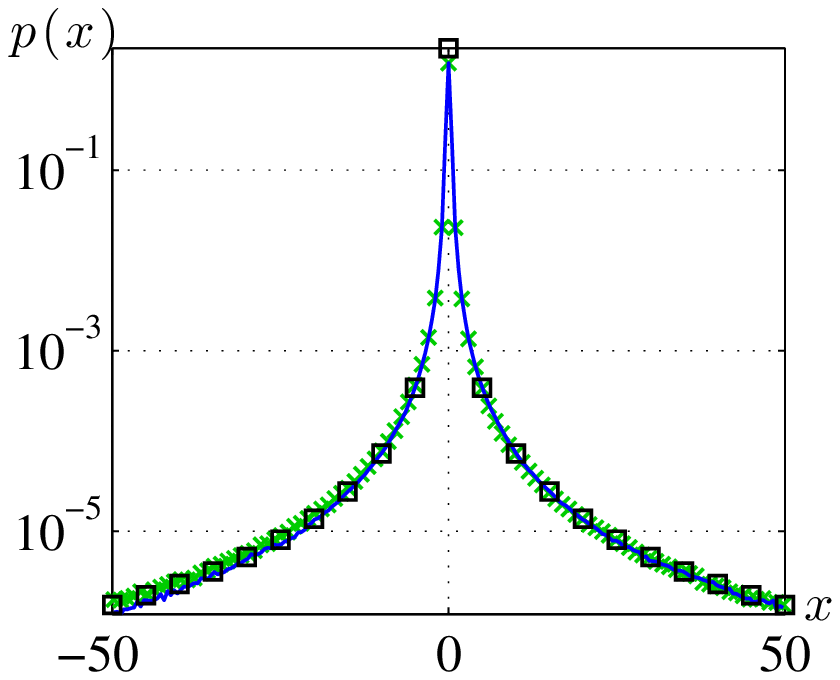}
\includegraphics[width=0.323\textwidth]{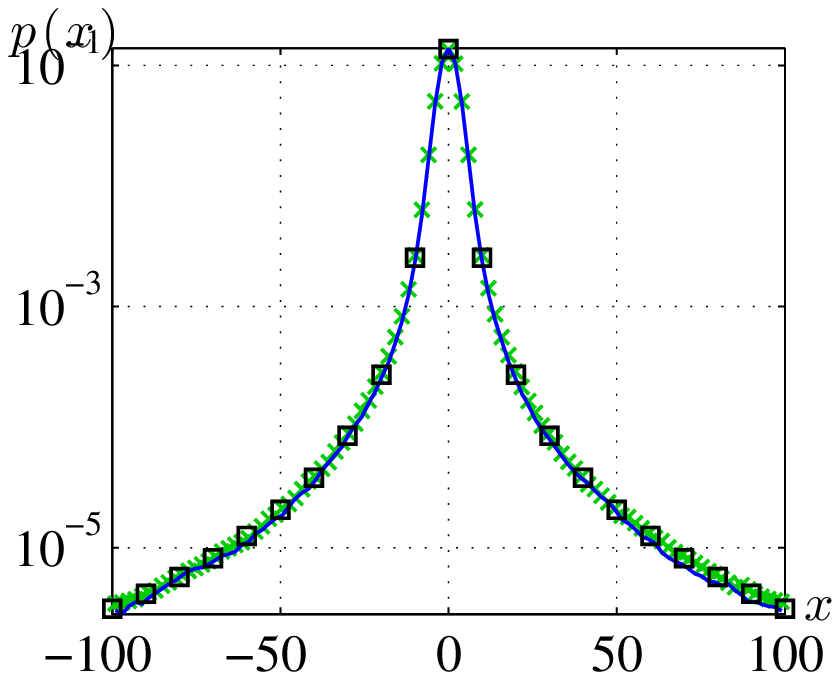}
\caption{Normalized histogram-based estimates of the PDF with symmetric forcing ($\beta = 0$) of $x_{t}$ from (\ref{linsys1}, \ref{linsys2}) (blue line) and the (L) approximation $\xi_{t}$ of (\ref{Lapprox_linear}) (green crosses). Left/Middle: As in Fig. \ref{linhisto_intro} with $\alpha = 1.9/1.4$. Right: As in Fig. \ref{linhisto_intro}, but with $a = 0.7,\,b = 2,\,\alpha = 1.7$ and $\epsilon = 0.1$. The numerically evaluated density using the MATLAB package STBL for $\xi$ in (\ref{Lapprox_linear}) is shown by black squares.} \label{linhisto}
\end{figure}

\begin{figure}[h]
\centering
\includegraphics[width=0.323\textwidth]{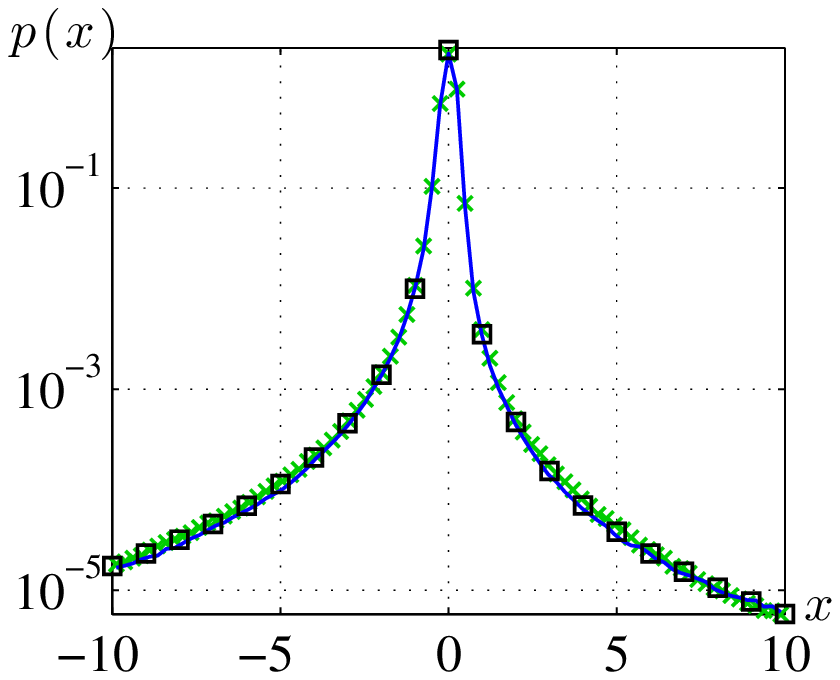}
\includegraphics[width=0.323\textwidth]{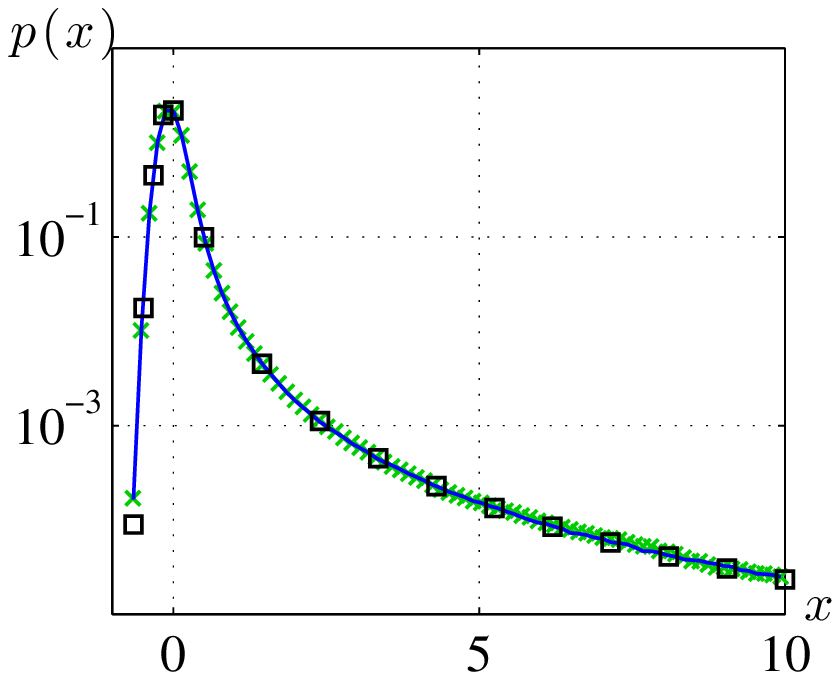}
\includegraphics[width=0.323\textwidth]{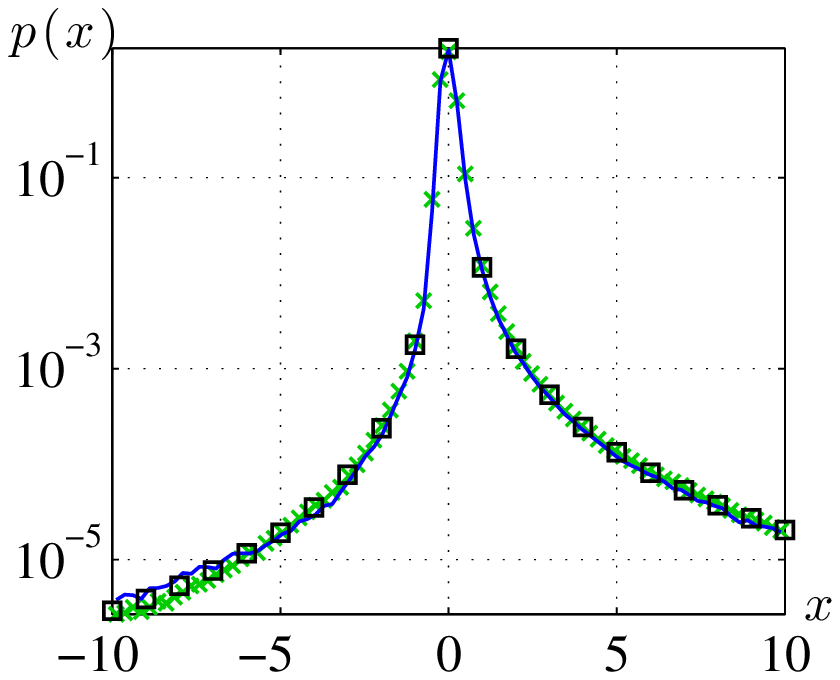}
\caption{Normalized histogram-based estimates of the PDF with $\beta \ne 0$ of $x_{t}$ from (\ref{linsys1}, \ref{linsys2}) (blue line). Also plotted is the histogram for $\xi_{t}$ of (\ref{Lapprox_linear}) (green crosses). For all figures: $a = 0.2,\,b = 0.7,\,c = 1,\,\alpha = 1.7$. Left/Middle/Right: $\epsilon = 0.01/0.01/0.1,\,\beta = -0.5/1/0.75$. The numerically evaluated density using the MATLAB package STBL for $\xi$ in (\ref{Lapprox_linear}) is shown in black squares.} \label{linhisto_skew}
\end{figure}

In Fig. \ref{linhisto}, we display the numerically simulated PDF of $x_{t}$ and $\xi_{t}$ with dynamics given by  (\ref{linsys1}, \ref{linsys2}) and (\ref{Lapprox_linear}) respectively for different values of $\alpha < 2$. We also consider cases with asymmetric noise (i.e. $\beta \ne 0$) in Fig. \ref{linhisto_skew}. Comparing the numerically generated densities of $\xi_{t}$ (\ref{Lapprox_linear}) with those of $x_{t}$ from (\ref{linsys1}, \ref{linsys2}) we see that the density obtained by the (L) approximation is almost indistinguishable from the density of the full system. In these figures, we also plot the numerically computed distribution as determined using the STBLPDF command in the STBL package for MATLAB \cite{STBL}. This code computes the PDF by numerically evaluating an integral representation of the $\alpha$-stable distribution as given in Theorem 1 of \cite{Nolan1997}. The densities estimated from numerically simulated time series are in excellent agreement with the numerically computed densities as determined from the integral representation.

We also compare the autocodifference functions (AFs) \cite{Taqqu1994} of the simulated trajectories for the full and averaged systems. This quantity is the generalization of the autocovariance function to variables that may not have a finite variance. In the case where $\alpha = 2$, the autocovariance function is a standard measure of the degree of linear dependence between time-lagged states of a time series. However, when $\alpha < 2$, the autocovariance function is not defined and so instead we must use the analogous AF to get a sense of this linear dependence. The AF for a continuous-time stochastic process $z_{t}$ is defined in terms of the CF of $z_{t}$, 
\begin{equation}
A_{z}(\tau) = \log\left[\expected{\exp(i(z_{t+\tau} - z_{t}))}\right] - \log\left[\expected{\exp(iz_{t+\tau})}\right] - \log\left[\expected{\exp(-iz_{t})}\right], \label{eqn:autocodiff}
\end{equation}
where $\tau$ is the lag time. In the case of $\xi_{t}$ given by (\ref{Lapprox_linear}), the AF can be explicity evaluated,
\begin{align}
{A}_{\xi}(\tau) &= \frac{(ab)^{\alpha}}{\alpha(1 - ac)}\left[1 + \exp(-\alpha(1 - ac)\tau) - |1 - \exp(-(1 - ac)\tau)|^{\alpha} \right] \label{eqn:autocodiff_linear} \\& \quad- i\beta\tan\left(\frac{\pi\alpha}{2}\right)\frac{(ab)^{\alpha}}{\alpha(1 - ac)}\left[ (1 - \exp(-\alpha(1 - ac)\tau)) - |1 - \exp(-(1 - ac)\tau)|^{\alpha}\right]. \nonumber
\end{align}
We note that when $\beta = 0$, the imaginary part of the AF is equal to zero and (\ref{eqn:autocodiff_linear}) reduces to the autcodifference function given in \cite{Taqqu1994}. If $\alpha = 2$, then $A_{\xi}$, is identical to the autocovariance function of $\xi$. The existence of a nonzero imaginary part of the AF for $\beta \ne 0$ is a consequence of definition (\ref{eqn:autocodiff}). It has no obvious interpretation in terms of the standard autocovariance function or the dynamics of the corresponding process.

\begin{figure}[h]
\centering
\includegraphics[width=0.323\textwidth]{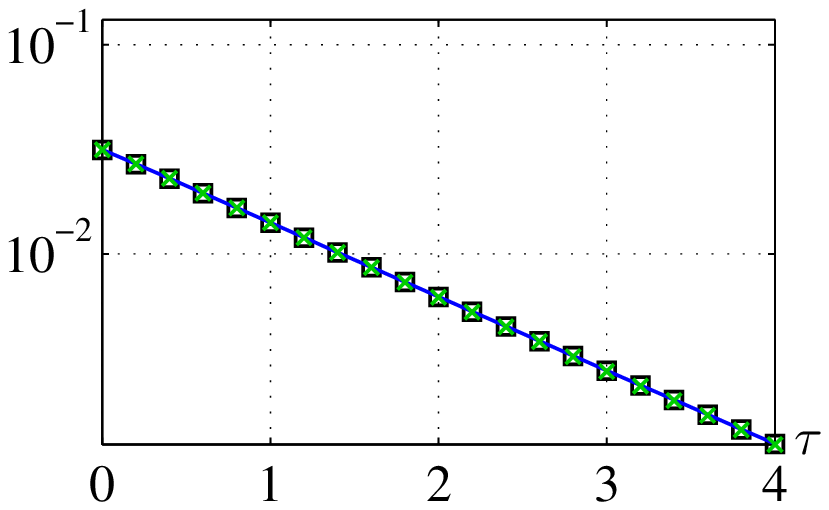}
\includegraphics[width=0.323\textwidth]{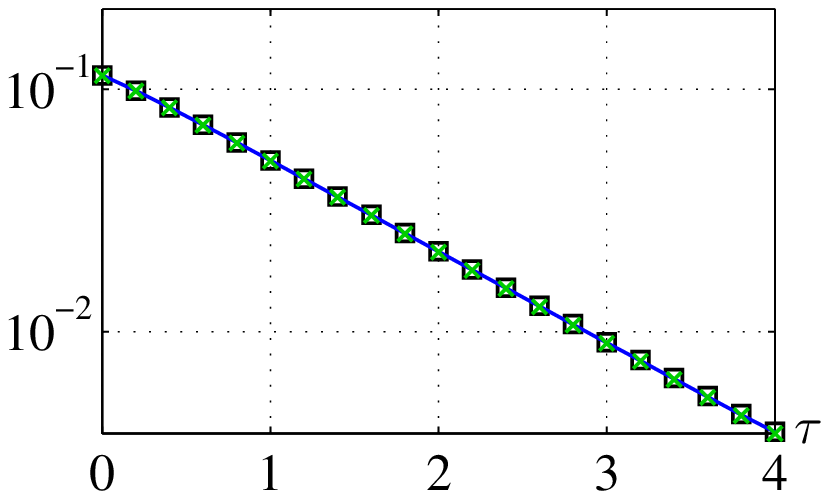}
\includegraphics[width=0.323\textwidth]{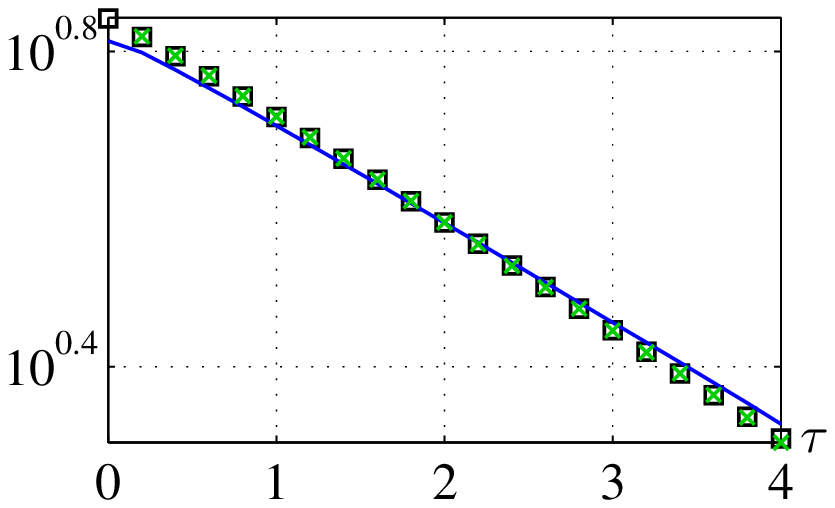}
\caption{Logarithmic plot of the numerically estimated AF for $x_{t}$ of (\ref{linsys1}, \ref{linsys2}) (blue line) and $\xi_{t}$ of (\ref{Lapprox_linear}) (green crosses) where $\beta = 0$. The analytically derived autocodifference $A_{\xi}(\tau)$ in (\ref{eqn:autocodiff_linear}) is also displayed (black squares). The parameters for the left, middle and right panels correspond to those from Fig. \ref{linhisto}.} \label{lin_autocodiff}
\end{figure}

We estimate the AF for the sample trajectories used to generate the density estimates given in Fig. \ref{linhisto} and offer a brief summary of the estimation method in Appendix \ref{app:numericalACD}. The results are displayed in Fig. \ref{lin_autocodiff}. We see that for $\epsilon \ll 1$, there is virtually no difference between the estimated AFs of the full and averaged systems, and both are comparable to the analytically derived result (\ref{eqn:autocodiff_linear}). For $\epsilon = 0.1$, small differences between the AFs of the full and averaged systems are apparent but not unexpected given the relatively large value of $\epsilon$ which is on the edge of the asymptotic regime.

We also consider the case of skewed noise (not shown) for $\epsilon \ll 1$ and see that the modulus of the estimated AFs of the full and reduced systems are indistinguishable for various values of $\beta$ and $\epsilon = 0.01$ with some slight differences between the full and averaged systems for $\epsilon = 0.1$. Similar agreement is observed for the real and imaginary components of the AF compared separately. From these results, we see that (\ref{Lapprox_linear}) is statistically similar to $x_{t}$ of (\ref{linsys1}, \ref{linsys2}) not only in stationary distribution, but also in autocodifference.

\subsection{Nonlinear system 1: Bilinear slow dynamics}
\label{subsec:nonlinsys1}
The first nonlinear system we study is a fast-slow system with a bilinear term in the slow equation:
\begin{align}
dx_{t} &= \left( c - x_{t} + \epsilon^{-\gamma}x_{t}y_{t}\right)\,dt, \quad x_{0} > 0, \label{nonlin1a}\\
dy_{t} &= -\frac{y_{t}}{\epsilon{a}}\,dt + \frac{b}{\epsilon^{1/\alpha}}\,dL_{t}^{(\alpha,0)}, \quad y_{0} = 0. \label{nonlin1b}
\end{align}
where $c,\,{a} > 0$. The bilinearity in the slow equation has the effect of restricting the dynamics of $x_{t}$ to positive values since the increments $dx_{t}/dt \rightarrow c > 0$ as $x \rightarrow 0^{+}$. We compute the (N+) approximation, ${X}_{t}$, using  (\ref{N+approx}),
\begin{equation}
d{X}_{t} = (c -{X}_{t})\,dt + {a}b{X}_{t}\diamond dL_{t}^{(\alpha,0)},\quad {X}_{0} = x_{0}. \label{nonlin1_N+}
\end{equation}
For the sake of comparison, we also give the (L) approximation using (\ref{Lapprox}),
\begin{equation}
x_{t} \approx \bar{x}_{t} + \xi_{t} \quad \mbox{where} \quad \begin{cases} d\xi_{t} &= -\xi_{t}\,dt + {a}b|\bar{x}_{t}|\,dL_{t}^{(\alpha,0)}, \quad \xi_{0} = x_{0} - c \\
\bar{x}_{t} &= c + (x_{0} - c)\exp(-t). \end{cases}
\label{nonlin1_L}\\
\end{equation}

\begin{figure}[h]
\centering
\includegraphics[width=0.323\textwidth]{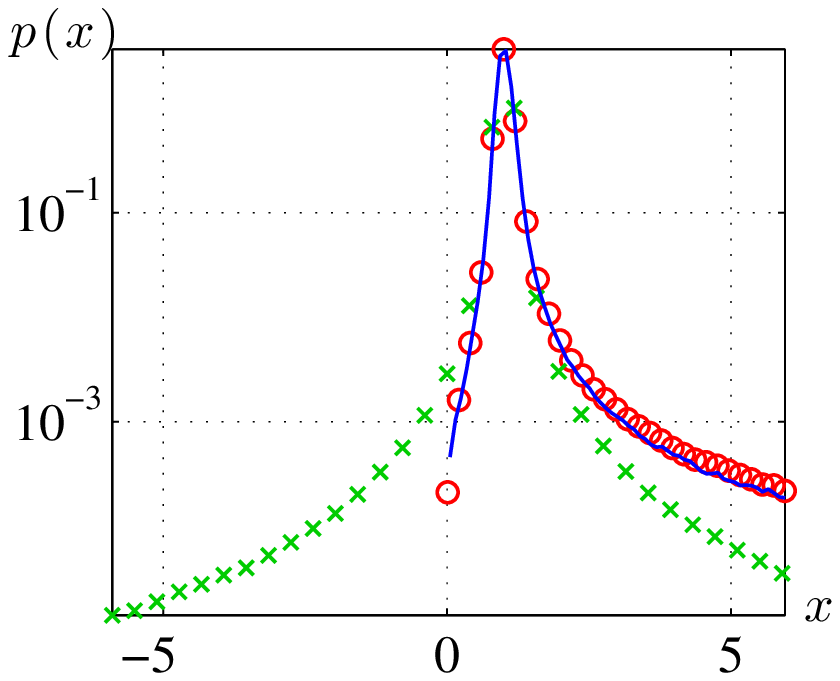}
\includegraphics[width=0.323\textwidth]{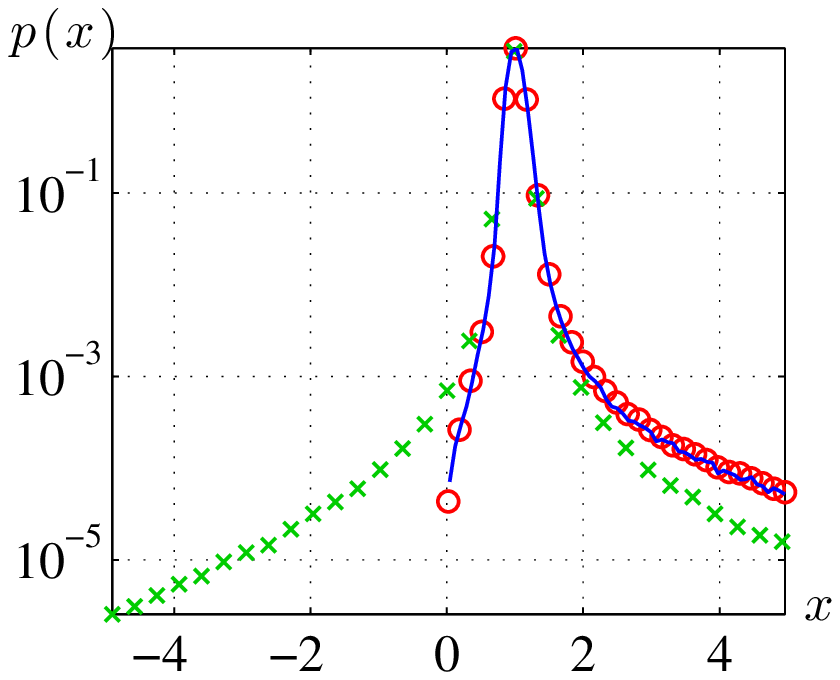}
\includegraphics[width=0.323\textwidth]{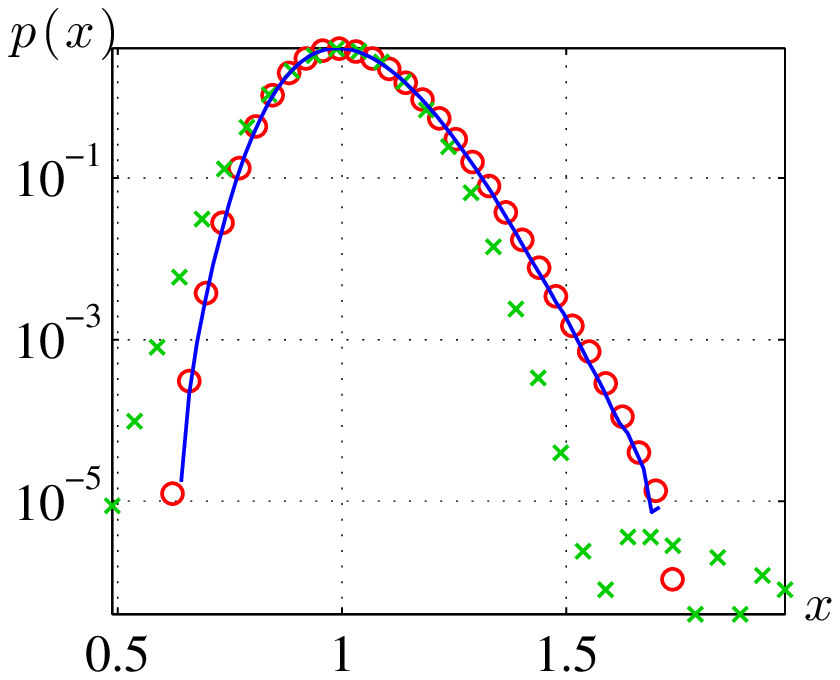}
\caption{Normalized histogram-based estimates of the PDF of $x_{t}$ for the full nonlinear system (\ref{nonlin1a}, \ref{nonlin1b}) (blue line), the (L) approximation $\bar{x}_{t} + \xi_{t}$ (\ref{nonlin1_L}) (green crosses) and the (N+) approximation ${X}_{t}$ (\ref{nonlin1_N+}) (red circles). For all figures: ${a} =1,\,b = 0.1,\,c = 1,\,\epsilon = 0.01$ and $x_{0} = c$. Left/Middle/Right: $\alpha = 1.7/1.9/2$.} \label{nonlinhisto1}
\end{figure}

We compare the densities of the full system (\ref{nonlin1a}, \ref{nonlin1b}) and the approximations (\ref{nonlin1_N+}) and (\ref{nonlin1_L}) via numerical simulation in Fig. \ref{nonlinhisto1} demonstrating the expected superiority of the (N+) approximation over the (L) approximation. The simulation methods are described in Appendix \ref{app:simulating_SDEs}, (It\={o} for (\ref{nonlin1a}, \ref{nonlin1b}) and (\ref{nonlin1_L}) and Marcus for (\ref{nonlin1_N+})). Particularly noteworthy is that the (N+) approximation, like the full system, does not cross the $x = 0$ boundary. Given the multiplicative nature of the noise of the (N+) approximation, it is clear that the positivity of the slow process is preserved in the approximation, as the dynamics of $X_{t}$ (3.9) satisfy $dX_{t}/dt \rightarrow c > 0$ for $X_{t} \rightarrow 0$. The dynamics of $x_{t}$ in the full system (\ref{nonlin1a}, \ref{nonlin1b}) also satisfy the same positivity constraint, as shown above. In contrast, since the (L) approximation is linear drift with additive noise, it is possible for the (L) approximation to cross zero. However, even though the (L) approximation does cross $x = 0$, it is a reasonable approximation of the full system near $\bar{x}_{t}$ (the mean of $x_{t}$), and thus may be suitable for local approximations. Although we have only presented results for $\beta = 0$, the approximations for $\beta \ne 0$ have similar quality (see Appendix \ref{app:simulation_Marcus} for a discussion on the numerical method).

\begin{figure}[h]
\centering
\includegraphics[width=0.323\textwidth]{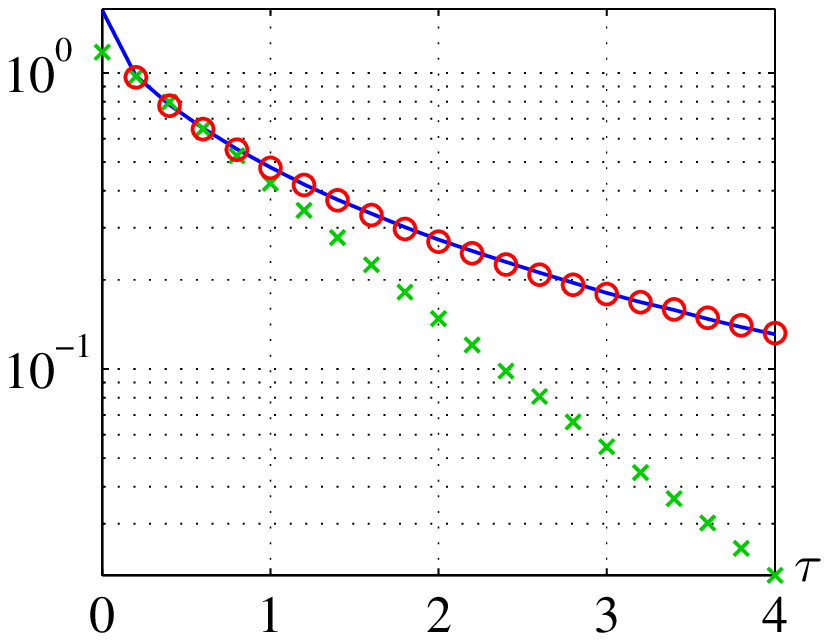}
\includegraphics[width=0.323\textwidth]{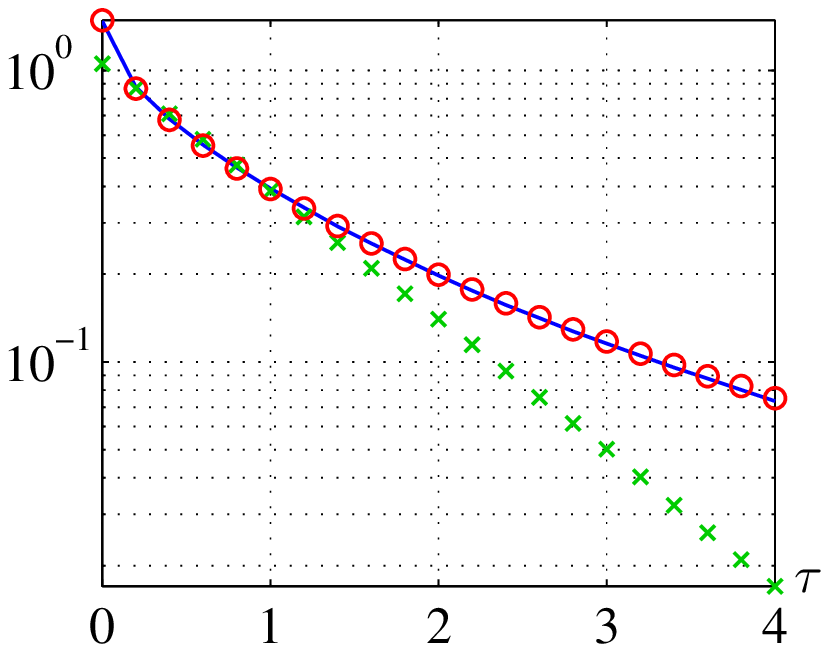}
\includegraphics[width=0.323\textwidth]{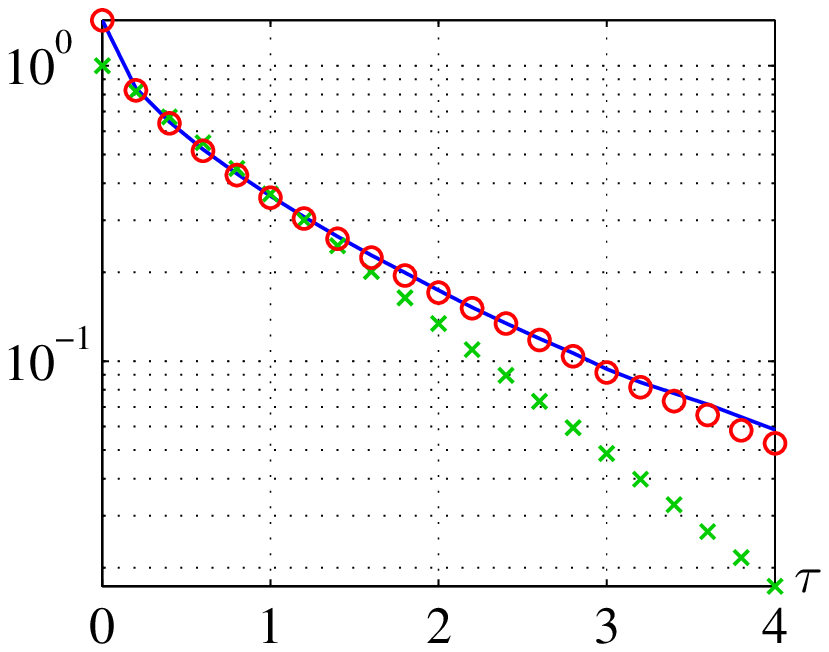}
\caption{Logarithmic plot  of the estimated AF for $x_{t}$ of the full nonlinear system (\ref{nonlin1a}, \ref{nonlin1b}) (blue line). Also plotted are the real parts of the autocodifferences of the (L) approximation, $\bar{x}_{t} + \xi_{t}$, (\ref{nonlin1_L}) (green crosses) and the (N+) approximation, ${X}_{t}$,  (\ref{nonlin1_N+}) (red circles). The parameters for the 3 panels correspond to those from Fig. \ref{nonlinhisto1}.} \label{nonlin1_autocodiff}
\end{figure}

Figure \ref{nonlin1_autocodiff} shows that the AF of the (L) approximation, which exhibits exponential decay, differs from that of the (N+) approximation and the full system. The slope of the AF for the (L) approximation for small time lags is shallower than that of (N+), indicating less rapid "memory loss" over short time intervals. However, for longer times, the (L) approximation appears to have more rapid memory loss than both $x_{t}$ and ${X}_{t}$ of the full and (N+) approximation systems. The autocodifference of the (N+) approximation is very similar to that of $x_{t}$ of the full system, indicating that ${X}_{t}$ is an excellent approximation to $x_{t}$ in both stationary PDF and AF.

For brevity we do not give figures of the AFs for the remaining examples as they give essentially the same result as those in Fig. \ref{nonlin1_autocodiff}; that the AFs of the (N+) approximation and the full system correspond very well for sufficiently small $\epsilon$, while the AF of the (L) approximation differs. 

\subsection{Nonlinear system 2: State-dependent reversion}
\label{subsec:nonlinsys2}
\begin{figure}[h]
\centering
\includegraphics[width=0.32\textwidth]{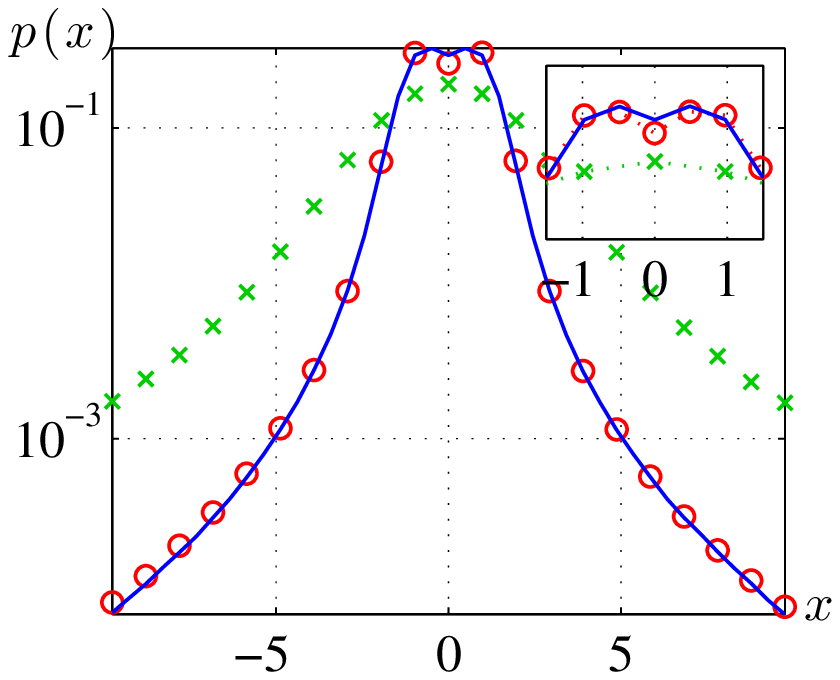}
\includegraphics[width=0.32\textwidth]{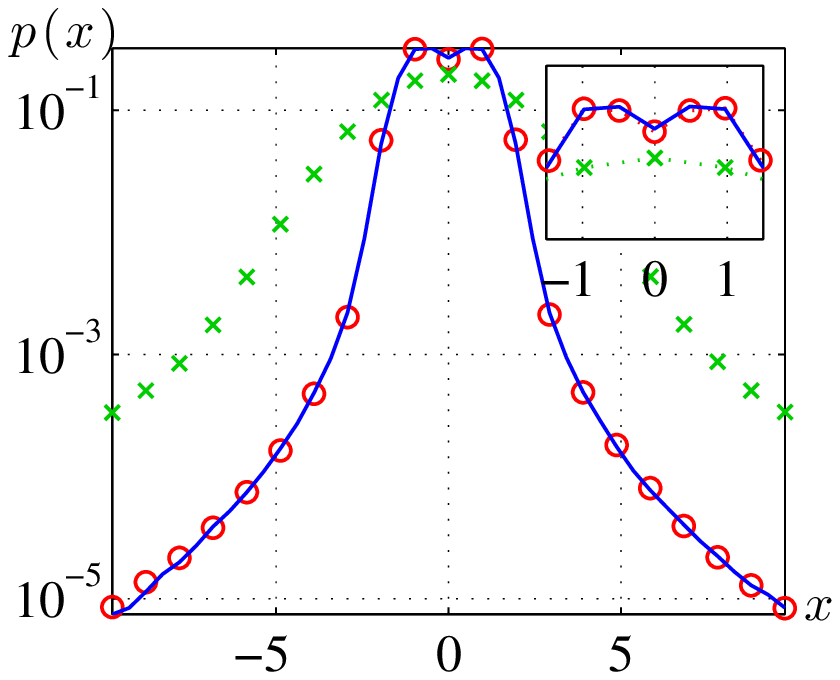}
\includegraphics[width=0.32\textwidth]{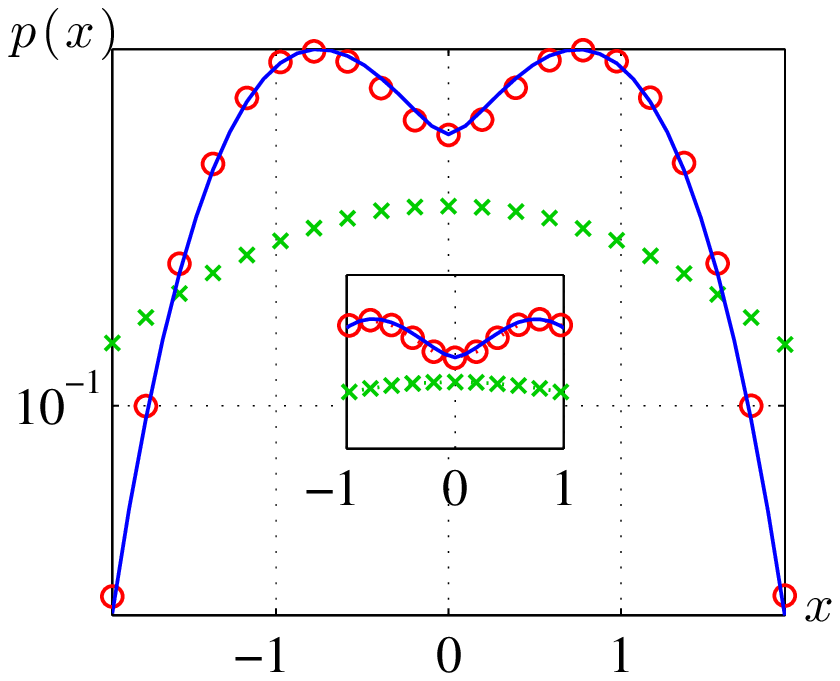}
\caption{Normalized histogram-based estimates of the PDF of $x_{t}$ for the full system (\ref{nonlin2a}, \ref{nonlin2b}) (blue line), the (L) approximation $\bar{x}_{t} + \xi_{t}$ (\ref{nonlin2_L}) (green crosses), and the (N+) approximation, ${X}_{t}$, (\ref{nonlin2_N+}) (red circles). For all plots, $b = 2$ and $x_{0} = 0$. Left: $\epsilon = 0.1,\,\alpha = 1.6$. Middle: $\epsilon = 0.01,\,\alpha = 1.9$. Right: $\epsilon = 0.01,\,\alpha = 2$. (Insets are on a linear scale)} \label{nonlinhisto2}
\end{figure}

For the second nonlinear system we consider, there is a nonlinear term in the equation for the fast variable $y_{t}$,
\begin{align}
dx_{t} &= \left( -x_{t} + \epsilon^{-\gamma}y_{t}\right)\,dt, \quad x_{0} \in \mathbb{R}, \label{nonlin2a}\\
dy_{t} &= -\frac{(1 + |x_{t}|)}{\epsilon}y_{t}\,dt + \frac{b}{\epsilon^{1/\alpha}}\,dL_{t}^{(\alpha,0)}, \quad y_{0} = 0. \label{nonlin2b}
\end{align}
For larger values of $|x|$, there is stronger reversion to zero of the fast variable $y_{t}$. Computing the (N+) approximation using (\ref{N+approx}) and comparing to the (L) approximation (\ref{Lapprox}) gives
\begin{align}
\mbox{(N+):}& \quad x_{t} \approx {X}_{t}  \quad \mbox{where} \quad d{X}_{t} = -{X}_{t}\,dt + \frac{b}{(1 + |{X}_{t}|)}\diamond dL_{t}^{(\alpha,0)}, \quad {X}_{0} = x_{0}, \label{nonlin2_N+}\\
\mbox{(L):}&\quad x_{t} \approx \bar{x}_{t} + \xi_{t} \quad \mbox{where} \quad \begin{cases} \bar{x}_{t} = x_{0}\exp(-t),\\  d\xi_{t} = -\xi_{t}\,dt + \frac{b}{(1 + |\bar{x}_{t}|)}\,dL_{t}^{(\alpha,0)}, \quad \xi_{0} = 0. \end{cases} \label{nonlin2_L}
\end{align}
Figure \ref{nonlinhisto2} shows the estimates of the stationary PDFs for (\ref{nonlin2a}, \ref{nonlin2b}) and the approximations (\ref{nonlin2_N+}) and (\ref{nonlin2_L}) obtained by simulation. Once again, we see that the (N+) approximation gives an excellent approximation to the density of the full system. By comparison, the (L) approximation is poor as it does not capture the bimodal nature of the full system in the given parameter regimes and overestimates the probability mass in the tails of the distribution. In failing to capture the bimodal nature of  (\ref{nonlin2a}, \ref{nonlin2b}), the (L) approximation does not perform well even as a local approximation of the full system near the mean, $x = 0$. The (N+) approximation accurately captures both the tails and the bimodal nature of the full system.

\subsection{Nonlinear system 3: Strong slow nonlinearity}
\label{subsec:nonlinsys3}
\begin{figure}[h]
\centering
\includegraphics[width=0.32\textwidth]{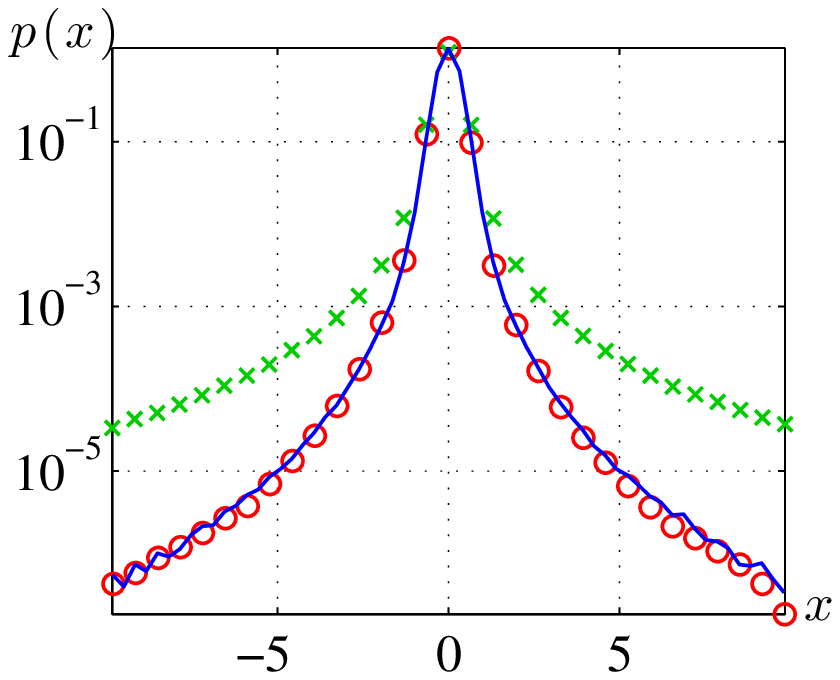}
\includegraphics[width=0.32\textwidth]{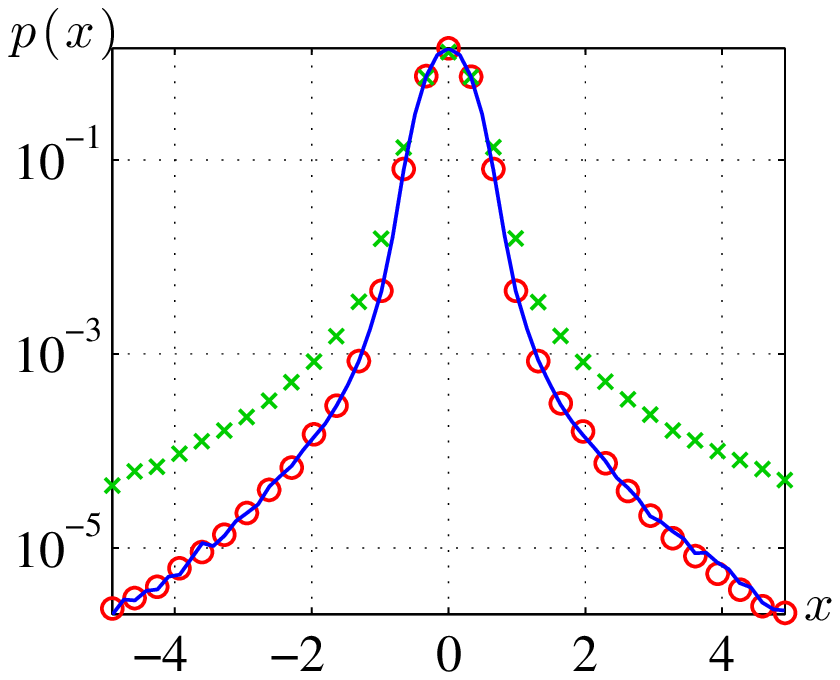}
\includegraphics[width=0.32\textwidth]{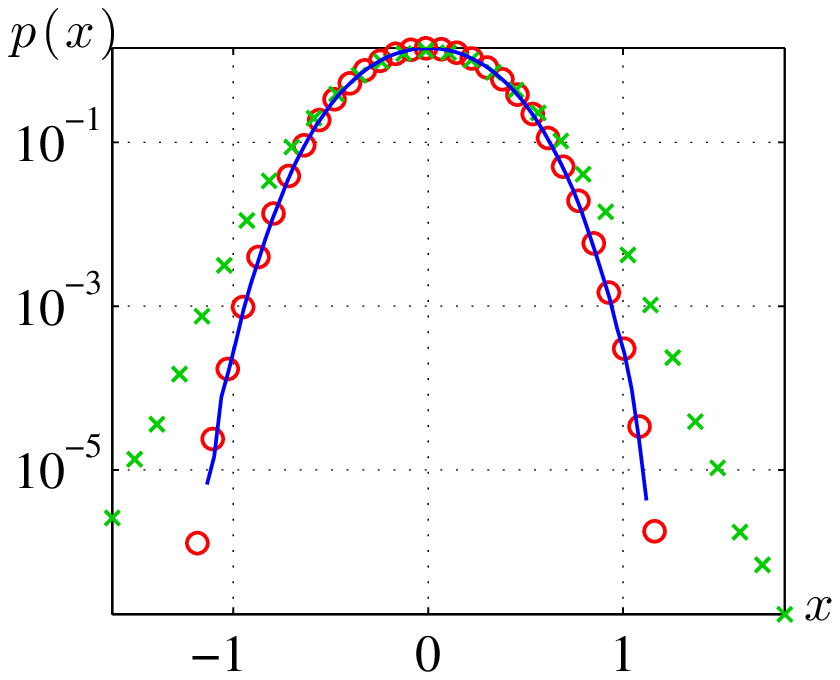}
\caption{Normalized histogram-based estimates of the PDF of the $x_{t}$ - variable for (\ref{nonlin3a}, \ref{nonlin3b}) (blue line), $\bar{x}_{t} + \xi_{t}$ for the system (\ref{nonlin3_L}) (green crosses), and ${X}_{t}$ for (\ref{nonlin3_N+}) (red circles). For all figures: $\sigma = 0.3,\,{a} =1,\,\epsilon = 0.01,\,x_{0} = 0$. Left/Middle/Right: $\alpha = 1.7/1.9/2.0$.}
\label{nonlinhisto3}
\end{figure}

In the third nonlinear model, the drift term of the slow variable contains a cubic term in $x_{t}$,
\begin{align}
dx_{t} &= \left( - x_{t} -  x_{t}^{3} + \epsilon^{-\gamma}y_{t}\right)\,dt, \quad x_{0} = 0, \label{nonlin3a}\\
dy_{t} &= -\frac{y_{t}}{\epsilon{a}}\,dt + \frac{b}{\epsilon^{1/\alpha}}\,dL_{t}^{(\alpha,0)}, \quad y_{0} = 0. \label{nonlin3b}
\end{align}
where ${a} > 0$. The fast variable which is an OULp, appears linearly in the slow equation. For $x_{t}$ near 0, the system can be approximated as a multivariate OULp as the cubic term in (\ref{nonlin3a}) is small. As before, we compute the (N+) approximation (\ref{N+approx}) and compare it to the (L) approximation (\ref{Lapprox}), 
\begin{align}
\mbox{(N+):}& \quad x_{t} \approx X_{t}  \quad \mbox{where} \quad d{X}_{t} = (-{X}_{t}-{X}_{t}^{3})\,dt + {a}b\,dL_{t}^{(\alpha,0)}, \quad X_{0} = x_{0}, \label{nonlin3_N+} \\
\mbox{(L):}&\quad x_{t} \approx \bar{x}_{t} + \xi_{t} \quad \mbox{where} \quad \begin{cases} \bar{x}_{t} = \frac{x_{0}}{\sqrt{x_{0}^{2}(e^{2t} - 1) + e^{2t}}}, \\ d\xi_{t} = -\xi_{t}\,dt + {a}b\,dL_{t}^{(\alpha,0)}, \quad \xi_{0} = 0. \end{cases} \label{nonlin3_L}
\end{align}
From the results in Fig. \ref{nonlinhisto3}, we see that again the (N+) approximation captures the full nonlinearity of the slow dynamics (\ref{nonlin3a}) while the (L) approximation effectively gives a linearization of the drift dynamics near the deterministic equilibrium of (\ref{nonlin3a}) when $y_{t}$ is at its mean value of 0. Both approximations do well near the peak of the distribution at $x = 0$, but the (L) approximation deviates from the full system for $|x|$ large, so that the OULp approximation is not appropriate. The (N+) approximation, by comparison, is once again a good approximation over the whole domain.

\section{Conclusions}
\label{sec:conc}

In this paper, we determine  the (L) and (N+) appproximations for fast-slow sytems forced with $\alpha$-stable noise. These are reduced systems providing weak approximations for the slow variables, analogous to the established (L) and (N+) approximation for Gaussian noise forcing. Previous results for Gaussian noise have been based on the analysis of autocovariance \cite{Arnold2003, Mona11}, while other analyses have used the characteristic functions obtained from the Fourier transform of the FKE for linear fast-slow systems \cite{Srokowski2011}. Here we consider nonlinear systems, and since the autocovariance function is not available for systems with $\alpha$-stable noise, we use the characteristic functions for appropriately transformed variables as the basis for the approximations. We restrict our attention to nonlinear models which are linear in the fast variable, as the Fourier analysis necessary for determining the CF is straightforward in this case.

For the (L) and (N+) approximations that we derive, we find that the drift and diffusion coefficients of the dynamics have the same form up to a power of $\epsilon$ as the coefficients obtained in the (L) and (N+) approximations in the Gaussian case. The difference is in the treatment of the noise; obviously, for systems with $\alpha$-stable noise, the (L) and (N+) approximations also involve $\alpha$-stable noise. In addition, the general result for the (N+) approximation requires a Marcus interpretation for the  $\alpha$-stable case, while in the Gaussian case the interpretation of the noise terms is Stratonovich. For $\alpha = 2$, corresponding to Gaussian noise with no jumps, the Marcus interpretation reduces to the Stratonovich, demonstrating consistency.

The derivation of the stochastic averaging methods presented here depends on the scaling relationships with regard to the timescale separation parameter $\epsilon$. In particular, the value of the exponent $\gamma$ in the dynamical equations of the unreduced system, 
\begin{equation*}
\begin{cases}
dx_{t} = f_{1}(x_{t})\,dt + \epsilon^{-\gamma}f_{2}(x_{t})y_{t}\,dt, \\ \epsilon \, dy_{t} = \epsilon^{\gamma}g_{1}(x_{t}) + g_{2}(x_{t})y_{t}\,dt + \epsilon^{\gamma}b\,dL_{t}^{(\alpha,\beta)}, 
\end{cases}
\end{equation*}
is critically important. If $\gamma < 1 - 1/\alpha$, the stochastic terms of the (L) and (N+) approximations, which have the coefficients proportional to $\epsilon^{\gamma - (1 - 1/\alpha)}$, diverge in the formal $\epsilon \to 0$ limit. If $\gamma > 1 - 1/\alpha$, the (L) and (N+) approximations are asymptotically equivalent to the (A) approximation as $\epsilon \rightarrow 0$, which does not include any stochastic effects. For non-zero $\epsilon = o(1)$ and $\gamma > 1-1/\alpha$, the (L) and (N+) approximations are asymptotic approximations in the weak sense for the dynamics of $x_{t}$. While we also have this approximation for $\gamma < 1 - 1/\alpha$, the dynamics of $x_{t}$ become dominated by the noise in this case and it becomes questionable that $x_{t}$ can be considered a slow process relative to $y_{t}$. For the critical value $\gamma = 1 - 1/\alpha$, the intensity of the noise in the (L) and (N+) approximations does not depend on $\epsilon$ and the stochastic corrections persist as the timescale separation parameter $\epsilon$ is made arbitrarily small.

We assumed that $g_{2}(x) < 0$ is $O(1)$ to ensure contracting dynamics so that the distribution of $y$ conditioned on $x$ reaches stationary behaviour for $y$ on a fast time scale. We expect that other systems with  $g_{2}(x) > 0$ for some $x$ also have this property, but we have not considered them here in order to concentrate on developing the approach for deriving the reduced system, rather than on confirming that conditionally stationary behaviour of $y$ is realized on a fast time scale. In addition, we expect that caution is required in the cases where the distribution of y conditioned on $x$ is multimodal  (e.g. multistability for $y$ in the absence of noise).


We compare numerical results for several fast-slow systems (both linear and nonlinear) and their corresponding (L) and (N+) approximations for the slow variables, computing both the long time probability density functions and the autocodifference functions (analogous to the autocovariance functions for $\alpha$-stable processes). We observe good agreement between the full system and the approximations with some expected limitations. The (L) approximation gives a good approximation for linear fast-slow systems as well as a local approximation near the mean when a linear approximation is valid there. The (N+) approximation captures the nonlinear behaviour of both the long time probability density functions and the autocodifference functions, showing good agreement with the full system over all values for the slow variable. However, the (N+) approximation is more difficult to simulate, primarily because of the Marcus stochastic terms. Considerations of numerical stability of nonlinear systems may also be an issue, requiring higher order methods for accuracy, as discussed in Appendix \ref{app:simulation_Ito}. The choice to use the (L) or (N+) approximation depends on the balance of accuracy, range of values where the approximation is needed, and the amount of analytical and computational work required.

\begin{figure}[h]
\centering
\includegraphics[width=0.323\textwidth]{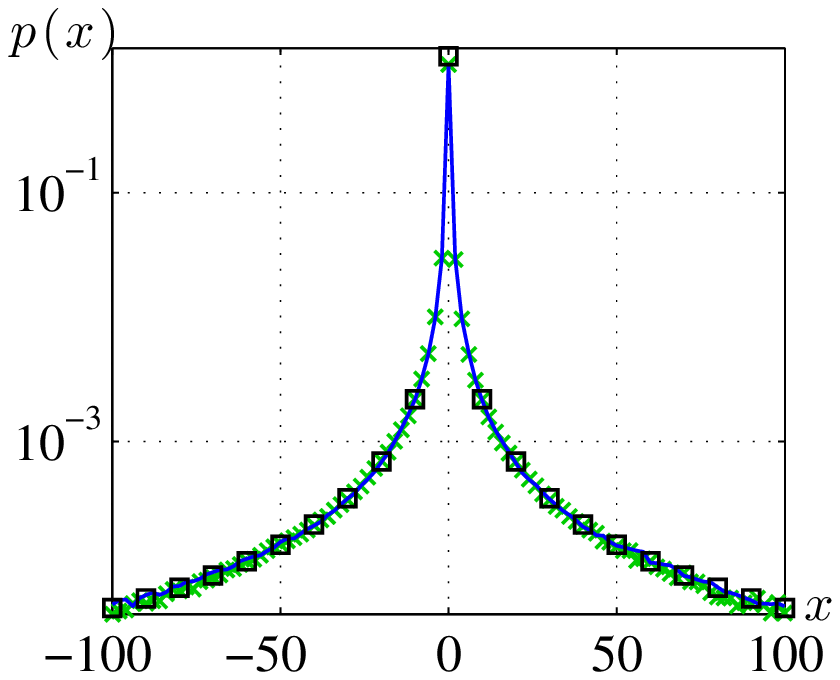}
\includegraphics[width=0.323\textwidth]{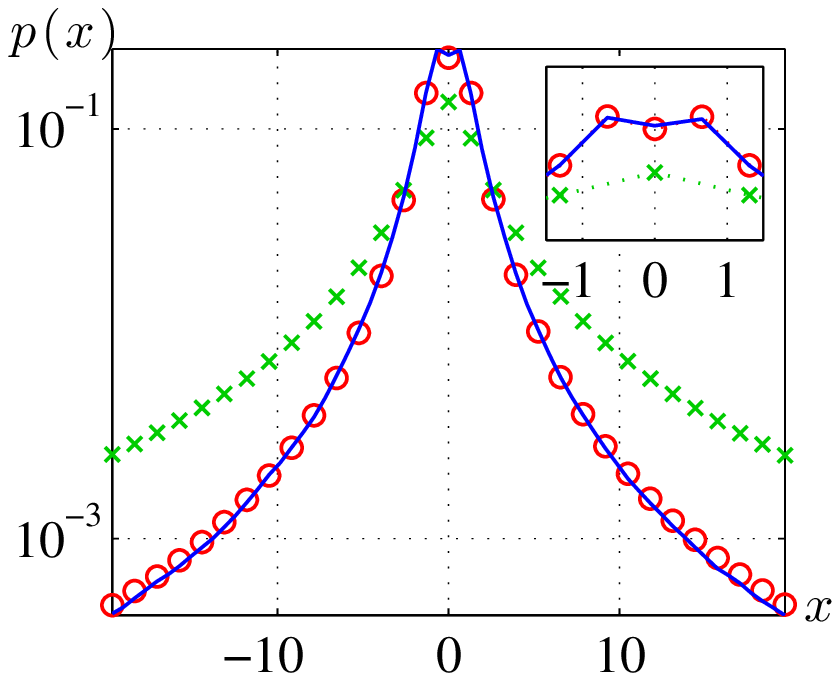}
\includegraphics[width=0.323\textwidth]{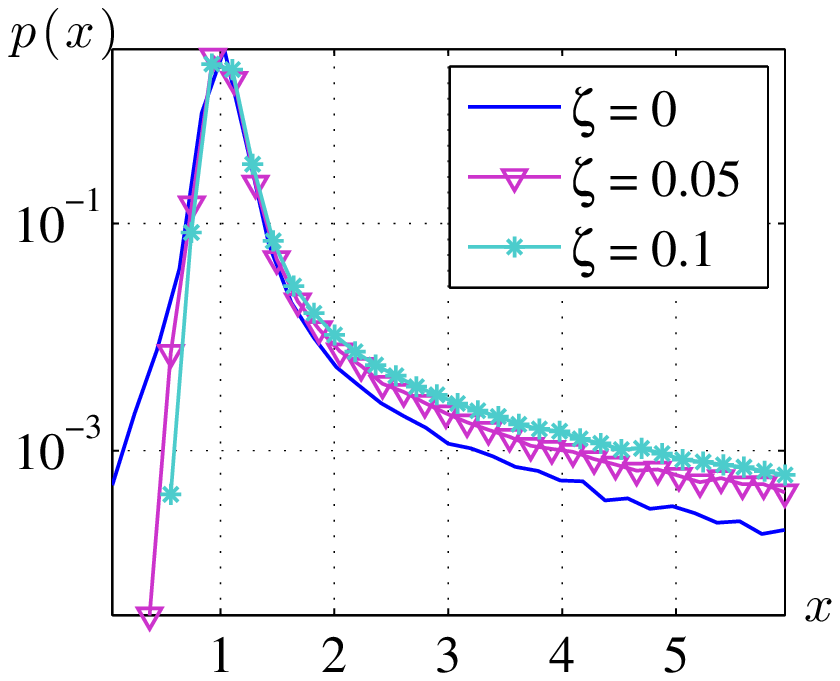}
\caption{Left \& Middle: Normalized histogram-based estimates of the PDF with symmetric $\alpha$-stable forcing for $\alpha < 1$, which is not explicitly covered in our analysis. (Left: As in Fig. \ref{linhisto_intro} with $(a,b,c,\epsilon,\alpha) = (0.2,0.7,1,10^{-2},0.7)$. Middle: As in Fig. \ref{nonlinhisto2} with $(\epsilon,\alpha) = (10^{-2}, 0.8)$). Right: Estimates of the PDF for $x_{t}$ where $dx_{t} = \left( c - x_{t} + \epsilon^{-\gamma}x_{t}y_{t}\left[1 + \zeta y_{t}\right]\right)dt$ and $y_{t}$ satisfies (\ref{nonlin1b}), $(a,b,c,\epsilon,\alpha) = (1,0.1,1,10^{-2}, 1.7)$. The values of $\zeta$ are indicated in the legend.}
\label{linhisto_intro_small_alpha}
\end{figure}

Our results suggest a number of valuable areas for future research. In our analysis, we required that $\alpha>1$ so that the mean of the fast process is defined.  We show numerical simulations in Fig. \ref{linhisto_intro_small_alpha} that indicate cases where the averaging approximations are reasonable for $\alpha < 1$. This suggests that the (L) and (N+) approximations in the case of $\alpha < 1$ would be of interest for future study. The analysis of that case would likely require the use of a location parameter instead of the mean. The development of a stochastic reduction process from alternative perspectives (e.g. operator theory) would require substantial work beyond the scope of this paper, given the lack of higher order moments for the processes considered here. A comparison with such approaches would be valuable, and we leave it for future work. Another area of extension would be the consideration of systems with more general nonlinearities, including nonlinearities in the fast variables. In that case an explicit analytical result may not be available via Fourier analysis, but we expect that the densities for some systems that are nonlinear in $y$ 
are similar in character to those from systems that are linear in $y$. In Fig. \ref{linhisto_intro_small_alpha}, right panel, we show an example that is quadratic in $y$, demonstrating that the numerical approximation for the density is similar in character to that of our nonlinear system 1 with bilinear form.  It may then be possible to provide a (semi)-analytical expression for the reduced system, but the nonlinearities would require methods beyond those used in the present study. Furthermore, in order to consider (N+) approximations for systems with asymmetric $\alpha$-stable noise where the Marcus map cannot be expressed in closed form, a numerical method needs to be developed. A computationally expensive approach would be simply to treat every stochastic increment as a jump and use the numerical method outlined in Appendix \ref{app:simulation_Marcus}, but ideally one would prefer to have a more efficient method. Extending the analysis to systems where the fast and slow subsystems are multi-dimensional, rather than scalar is another natural direction for our research, although this may not be a trivial endeavour. A detailed comparison with strong pathwise approximations, such as \cite{Xu2011}, could provide an interesting study in which to explore weak convergence rates and the connection between the two distinct types of averaging.

\appendix
\section{Table of notation}
The following table gives a short list of the most commonly used variables and symbols in this study.
{\footnotesize
\begin{center}
\begin{tabular}{c|l}
Symbol(s) & Representative of \\ \hline
$x, y$ & the slow, fast dynamical variables \\
$\bar{x}$ & the process satisfying the (A) approximation to $x$ (\ref{Aapprox}) \\
$\xi$ & the perturbation from $\bar{x}$ in the (L) approximation to $x$ (\ref{Lapprox}) \\
$X$ & the (N+) approximation to $x$ (\ref{N+approx}). \\
$v,\,V$ & continuous processes approximated asymptotically with a white noise \\ & (defined in (\ref{v_dynamics}), (\ref{eqn:v_dynamics_N+1}))\\
$p(y,v,t),q(h,V,t)$ & the marginal probability distribution functions of $(y_{t},v_{t})$, $(h_{t},V_{t})$ resp. \\
$\mathcal{F}[f]$ & the Fourier transform of the function $f$ (see (\ref{char_func}) and (\ref{char_func2})) \\
$l,m,\omega$ & the Fourier variables associated with $y, v$ (and $V$), $h$ resp.\\
$\psi,\,\phi$ & the characteristic functions (CFs) of the $(y,v),\,(h,V)$ system resp. \\
$\mu$ & a variable that weakly approximates $\mathcal{T}(x)$ where $\mathcal{T}$ satisfies (\ref{eqn:Ttransformation})\\
\end{tabular}
\end{center}}

\section{Evaluation of the characteristic function $\psi$}
\label{app:calculations}
In this appendix we show the explicit calculations involving the asymptotic evaluation of the integral (\ref{integral_result}) that contributes to (\ref{psi_integral}), 
\begin{equation}
\int_{0}^{{t}}|\Lambda(r)|^{\alpha}\Xi(\Lambda(r);\alpha,\beta)\,dr \label{app_integral}
\end{equation}
where $\Xi,\, \Lambda$ are defined in (\ref{Xi}), (\ref{Lambda}) respectively. We focus on two different temporal regimes, $t = O(\epsilon)$ and $t = O(1)$. 

In the case where $t = O(\epsilon)$, we make a change of variable $\epsilon y = r$ in (\ref{app_integral}),
\begin{equation}
\epsilon\int_{0}^{\tilde{t}}|\Lambda(\epsilon y)|^{\alpha}\Xi(\Lambda(\epsilon y);\alpha,\beta)\,dy, \quad \tilde{t} = t/\epsilon = O(1). \label{app_integral_small_t}
\end{equation}
We use a Taylor series approximation for $|\Lambda(\epsilon y)|^{\alpha}$ to expand the integrand in powers of $\epsilon$ for $0 < \epsilon \ll 1$,
\begin{align}
|\Lambda(\epsilon y)|^{\alpha} &= \left|l\exp\left(g_{2}y\right) - m\epsilon^{1-\gamma}\frac{f_{2}}{g_{2}}\left(1 - \exp\left(g_{2}y\right)\right)\right|^{\alpha} \\
&\sim |l|^{\alpha}\exp(\alpha g_{2}y) - \epsilon^{1-\gamma} \alpha |l|^{\alpha-1}\sgn{l} m\frac{f_{2}}{g_{2}}\left[\exp((\alpha - 1)g_{2}y) - \exp\left(\alpha g_{2}y\right)\right].
\label{Lambda_small_t}
\end{align}
The function $\Xi(s;\alpha,\beta) = 1 + i\beta\sgn{s}\tan(\pi\alpha/2)$ depends on $s$ only though $\operatorname{sgn}(s)$ and is therefore constant except when $s$ changes sign. In the $t = O(\epsilon)$ limit, $\Lambda(\epsilon y) =  l\exp(g_{2} y) + O(\epsilon^{1-\gamma})$ and so the value of $\Xi$ is determined by the sign of $l$ for $y_{t}$ and $\epsilon$ sufficiently small:
\begin{equation}
\Xi(\Lambda(\epsilon y);\alpha,\beta) = \Xi(l;\alpha,\beta). \label{Xi_small_t}
\end{equation}
Substituting (\ref{Lambda_small_t}), (\ref{Xi_small_t}) into (\ref{app_integral_small_t}), we obtain an approximation for (\ref{app_integral}) for $t = O(\epsilon)$,
\begin{align}
\int_{0}^{{t}}|\Lambda(r)|^{\alpha}&\Xi(r;\alpha,\beta)\,dr \approx -\frac{\epsilon |l|^{\alpha}\Xi(l;\alpha,\beta)}{\alpha g_{2}}\left(1 - \exp\left(g_{2}t/\epsilon\right)\right) \\ &+ \epsilon^{2-\gamma} |l|^{\alpha-1}\Xi(l;\alpha,\beta)\sgn{l} m\frac{f_{2}}{g_{2}^{2}}\left(\frac{\alpha}{\alpha - 1}[1 - \exp((\alpha - 1)g_{2}t/\epsilon)] - [1 - \exp\left(\alpha g_{2}t/\epsilon\right)]\right) + O(\epsilon^{3-2\gamma}).
\end{align}
Keeping only the leading order term in $l$, we get the approximation to (\ref{app_integral}) valid for $t = O(\epsilon)$ that is written in (\ref{integral_result}).

For $t = O(1)$, (\ref{app_integral}) depends on both global contributions from terms that are significant over the entire domain of integration, and on local contributions that are concentrated on an interval of length $O(\epsilon)$ near $t = 0$. Methods for handling such integrals can be found in \cite{Hinch1991}. 
For the global contributions, we neglect terms in $\Lambda(r)$ that are exponentially small for $0 < \epsilon \ll 1$ and $r = O(1)$ to get the leading order behaviours of $|\Lambda(r)|^{\alpha}$ and $\Xi(\Lambda(r);\alpha,\beta),$
\begin{equation}
\begin{cases} |\Lambda(r)|^{\alpha} \sim \left|-\epsilon^{1-\gamma}m\frac{f_{2}}{g_{2}}\right|^{\alpha} = \epsilon\left|\frac{f_{2}}{g_{2}}\right|^{\alpha}|m|^{\alpha}, \\
\Xi(\Lambda(r);\alpha,\beta) = \Xi(f_{2}m;\alpha,\beta) = \Xi(m;\alpha,\beta^{*}), 
\end{cases}
\label{glob_approx}
\end{equation}
where $\beta^{*}$ is defined in (\ref{beta_star}). Then we evaluate (\ref{app_integral}) to get the global contribution plus a remainder, $\mathcal{R}$.
\begin{equation}
\int_{0}^{{t}}|\Lambda(r)|^{\alpha}\Xi(\Lambda(r);\alpha,\beta)\,dr = \epsilon\left|\frac{f_{2}}{g_{2}}\right|^{\alpha}|m|^{\alpha}\Xi(m;\alpha,\beta^{*})t + \mathcal{R}. \label{integral_w_remainder}
\end{equation}
Then $\mathcal{R}$ is primarily the local contribution near $r = 0$ and can be written as
\begin{align}
\mathcal{R} &= \int_{0}^{{t}}|\Lambda(r)|^{\alpha}\Xi(\Lambda(r);\alpha,\beta)\,dr - \epsilon\left|\frac{f_{2}}{g_{2}}\right|^{\alpha}|m|^{\alpha}\Xi(m;\alpha,\beta^{*})t = \int_{0}^{t}\exp\left(\frac{\alpha g_{2}r}{\epsilon}\right)\Upsilon(r)\,dr \label{remainder} \\
\Upsilon(r) &=  \left|l - m\epsilon^{1-\gamma}\frac{f_{2}}{g_{2}}\left(\exp\left(\frac{g_{2}r}{\epsilon}\right) - 1\right)\right|^{\alpha}\Xi(\Lambda(r);\alpha,\beta)\\ & \quad - \epsilon\exp\left(-\frac{\alpha g_{2}r}{\epsilon}\right)\left|\frac{f_{2}}{g_{2}}\right|^{\alpha}|m|^{\alpha}\Xi(m;\alpha,\beta^{*}). \nonumber
\end{align}
The remainder is approximated by evaluating the integral asymptotically using Watson's Lemma \cite{Hinch1991} and keeping the first three terms,
\begin{align}
\mathcal{R} \sim -\frac{\epsilon|l|^{\alpha}}{\alpha g_{2}}&\Xi(l;\alpha,\beta) + \frac{\epsilon^{2-\gamma}f_{2}}{g_{2}^{2}}|l|^{\alpha-1}\Xi(l;\alpha,\beta)\sgn{l}m\left(\frac{1}{\alpha} + \frac{1}{\alpha^{2}} \right) \nonumber \\ & + \frac{\epsilon^{3-2\gamma}f_{2}^{2}}{g_{2}^{3}}|l|^{\alpha-2}m^{2}\Xi(l;\alpha,\beta)\left(\frac{1}{\alpha^{2}}- \frac{1}{\alpha}\right) - \frac{3\epsilon^{2}}{\alpha g_{2}}\left|\frac{f_{2}}{g_{2}}\right|^{\alpha}|m|^{\alpha}\Xi(m;\alpha,\beta^{*}). \label{watsons_eval}
\end{align}
Substituting (\ref{watsons_eval}) into (\ref{integral_w_remainder}), we obtain the approximation to (\ref{app_integral}) valid for  $t = O(1)$, given in (\ref{integral_result}).

\section{Numerical Methods}
\label{app:numerical}
\subsection{Simulating the SDEs}
\label{app:simulating_SDEs}
In this appendix, we describe the numerical methods used to simulate trajectories of the systems given in \S \ref{sec:comp_results}. Depending on whether the stochastic dynamics of a system are interpreted in the sense of It\={o} or Marcus, the stochastic forcing term is treated differently and give details below in Appendices \ref{app:simulation_Ito} and \ref{app:simulation_Marcus}.

For the systems studied in \S \ref{sec:comp_results}, we simulate the full systems of the form (\ref{canon1}, \ref{canon2}) using a time step of size, $\delta t$ (referred to as the simulation time step). Each simulated time series was sampled on a coarser time grid with sampling time step $Dt = 10^{-2}$. The simulation time step, $\delta t \leq \frac{1}{10}Dt$, with $Dt$ no longer than the characteristic time scale of the fast process. The choice of time steps was taken to ensure that the simulation time scale and the characteristic time scale of the fast process were well-separated and that we obtain a sufficiently long time series while keeping the number of data points to a ``computationally tractable'' number. Different simulation time steps were used for different systems to compromise between issues of numerical stability and computation time, but all systems, including the approximations, have the same sampling time step, $Dt$. 

For the (L) and (N+) approximations in \S \ref{sec:comp_results}, we used $\delta t = D t$. This choice samples the dynamics sufficiently well since the characteristic time scale of the processes is $O(1)$.  Thus we ensure that both the full and reduced systems are simulated over the same time intervals with the same sample data sizes. Table \ref{tab:simulation_params} gives the parameters used in our numerical simulations.

\begin{table}
\centering
\footnotesize
\begin{tabular}{l|c|l|l}
System & variable & characteristic time scale & $\delta t$ \\
\hline \hline
(\ref{linsys1}, \ref{linsys2}) & $y_{t}$ &$\epsilon \ge Dt$ & $10^{-3}$ \\
(\ref{Lapprox_linear}) & $\xi_{t}$ &$(1 - ac)^{-1} \ge 5/4$ & $Dt$ \\
\hline
(\ref{nonlin1a}, \ref{nonlin1b})$^{*}$ & $y_{t}$ &$\epsilon {a} \ge Dt$ & $2 \times 10^{-4}$\\
(\ref{nonlin1_N+})$^{**}$ & $X_{t}$ &$1$ & $Dt$\\
(\ref{nonlin1_L}) & $\xi_{t}$ &$1$ & $Dt$\\
\hline
(\ref{nonlin2a}, \ref{nonlin2b}) & $y_{t}$ &$\epsilon(1 + |x_{t}|)^{-1}$ ($\epsilon \ge Dt$) & $5 \times 10^{-4}$\\
(\ref{nonlin2_N+})$^{***}$ & $X_{t}$ &1 & $Dt$\\
(\ref{nonlin2_L}) & $\xi_{t}$ &1 & $Dt$\\
\hline
(\ref{nonlin3a}, \ref{nonlin3b}) & $y_{t}$ &$\epsilon {a} \ge Dt$ & $10^{-4}$\\
(\ref{nonlin3_N+})$^{*}$ & $X_{t}$ &$1$ & $Dt$ \\
(\ref{nonlin3_L}) & $\xi_{t}$ &$1$ & $Dt$\\
\end{tabular}
\caption{Simulation parameters for the examples considered in \S \ref{sec:comp_results}. The indicated variable determines the smallest characteristic time scale of the system, which is also indicated. All systems were sampled on the time step $Dt = 10^{-2}$. Single asterisks indicate the system was simulated using a higher-order method as described by (\ref{sim_scheme_alt}). Double/Triple asterisks indicate numerical Marcus integration without/with numerical evaluation of the Marcus map, $\theta$ (\ref{num_method_marcus}).}
\label{tab:simulation_params}
\end{table}

\subsubsection{Numerical simulations for \S \ref{sec:comp_results}:  It\={o} interpretation}
\label{app:simulation_Ito}
For all of the (L) approximations simulated in \S \ref{sec:comp_results}, $\xi_{t}$ satisfies an SDE of the form
\begin{equation}
d\xi_{t} = H(\xi_{t})\,dt + \kappa(\bar{x}_{t})\,dL_{t}^{(\alpha,\beta)}, \quad \xi_{0} = 0,
\label{eqn:Ito_sde_app}
\end{equation}
where $\bar{x}_{t}$ is treated as constant with respect to $\xi_{t}$, so that the noise term has the It\={o} interpretation. As we are interested in the stationary distribution for $\bar{x}_{t} + \xi_{t}$, we take the initial value for $\bar{x}_{t}$ to be its stationary value in all numerical simulations of the (L) approximation, and hence there is no need to numerically simulate $\bar{x}_{t}$ as it is constant. We denote the constant value $\kappa({\bar{x}}_{t})$ as $\bar{\kappa}$. To numerically simulate (\ref{eqn:Ito_sde_app}), in most cases we use an Euler-type numerical approximation,
\begin{equation}
\hat{\xi}_{n+1} = \hat{\xi}_{n} + H(\hat{\xi}_{n})\,\delta t + \bar{\kappa}\,\delta L_{n}^{(\alpha,\beta)}, \quad \hat{\xi}_{0} = 0, \quad n = 0,1,2,\dots \label{sim_scheme}
\end{equation}
where $\hat{\xi}_{n}$ denotes the numerical approximation to $\xi_{t}$ at time $t = n\cdot \delta t$ and the terms $\delta L_{n}^{(\alpha,\beta)} \sim \mathcal{S}_{\alpha}\left(\beta,(\delta t)^{1/\alpha}\right)$. Numerical schemes analogous to (\ref{sim_scheme}) are used for all systems in this paper, unless otherwise indicated in Table \ref{tab:simulation_params}. For two systems, the weak order 1.0 predictor-corrector method \cite{Kloeden1992} was required for additional accuracy: i) nonlinear system 1 (\ref{nonlin1a}, \ref{nonlin1b}) to ensure the $x = 0$ boundary was not crossed and ii) the (N+) approximation of system 3 (\ref{nonlin3_N+}) as the cubic term caused numerical trajectories to diverge when the standard Forward Euler step method was used. The predictor-corrector numerical scheme is given by
\begin{equation}
\begin{cases}
\hat{\xi}_{n+1} = \hat{\xi}_{n} + \frac{1}{2}\left(H(\hat{\xi}^\dag) + H(\hat{\xi}_{n})\right)\,\delta t +\bar{\kappa}\,\delta L_{n}^{(\alpha,\beta)}\\ \hat{\xi}^\dag = \hat{\xi}_{n} + H(\hat{\xi}_{n})\,\delta t + \bar{\kappa}\,\delta L_{n}^{(\alpha,\beta)} 
\end{cases},
\quad \hat{\xi}_{0} = 0.\label{sim_scheme_alt}
\end{equation}
In all cases the increments of $\alpha$-stable noise are given by $\delta L_{n}^{(\alpha,\beta)}$, with stability parameter $\alpha \in (0,1) \cup (1,2)$, skewness parameter $\beta$ and scale parameter $\sigma = (\delta t)^{1/\alpha}$. We use the algorithm in \cite{Dybiec2006} to generate these quantities as follows,
\begin{equation}
\delta L_{n}^{(\alpha,\beta)} = B_{\alpha,\beta}\frac{\sin(\alpha(\zeta_{1} + C_{\alpha,\beta}))}{\cos(\zeta_{1})^{1/\alpha}}\left( \frac{\cos(\zeta_{1} - \alpha(\zeta_{1} + C_{\alpha,\beta}))}{\zeta_{2}}\right)^{(1 - \alpha)/\alpha} \label{JW_LevyRV}
\end{equation}
where
\begin{equation}
\begin{cases}
B_{\alpha,\beta} &= (\delta t)^{1/\alpha}\left(\cos\left(\arctan\left(\beta \tan\left(\frac{\pi\alpha}{2}\right)\right) \right) \right)^{-1/\alpha}, \\ C_{\alpha,\beta} &= \frac{1}{\alpha}\arctan\left(\beta\tan\left(\frac{\pi\alpha}{2}\right) \right)
\end{cases}
\end{equation}
and $\zeta_{1},\zeta_{2}$ are respectively a uniform random variable on the open interval $(-\pi/2,\pi/2)$ and an exponential random variable with parameter 1. In the case $\alpha = 2$, the increments $\delta L_{n}^{(2,\beta)} = \delta L_{n}^{(2,0)}$ are Gaussian random variables with variance $2\delta t$ and the simulation scheme (\ref{sim_scheme}) is then an Euler-Maruyama method \cite{Kloeden1992}.

\subsubsection{Numerical methods for evaluating Marcus integrals}
\label{app:simulation_Marcus}

We summarize results of \cite{Asmussen2007book, Chechkin2014, Grigoriu2009, Li2013, Li2013Erratum, Sun2013, Protter1990} relevant for processes $z_{t}$ satisfying SDEs with multiplicative noise terms in the sense of Marcus which are indicated by the ``$\diamond$'' symbol, such as (\ref{nonlin1_N+}) and (\ref{nonlin2_N+}). These systems are of the form,
\begin{equation}
dz_{t} = H(z_{t})\,dt + \kappa(z_{t})\diamond dL_{t}^{(\alpha,\beta)}, \quad z_{0} \in \mathbb{R} \label{eqn:Marcus_sde_app}
\end{equation}
and can be written in integral form as
\begin{equation}
z_{t} = z_{0} + \int_{0}^{t}H(z_{s})\,ds + \int_{0}^{t}\kappa(z_{s})\diamond dL_{s}^{(\alpha,\beta)}. \label{eqn:Marcus_sde_app_inte}
\end{equation}
The first integral on the right-hand side of (\ref{eqn:Marcus_sde_app_inte}) corresponds to the drift terms and is numerically approximated using standard methods such as Euler, explicit trapezoidal methods.  Here $L_{t}^{(\alpha, \beta)}$ is a pure jump process and so the Marcus increments are expressed in terms of the Marcus map, $\theta(r; \Delta L_{s}, z_{s-})$ as follows,
\begin{align}
\int_{0}^{t}\kappa(z_{s})\diamond dL_{s}^{(\alpha,\beta)} = \sum_{s \le t}\left[\theta(1;\Delta L_{s},z_{s-}) - z_{s-} \right].\label{marcus_discret2}
\end{align}
where $\theta(r; \Delta L_{s}, z_{s-})$ satisfies
\begin{equation}
\frac{d\theta(r; \Delta L_{s}, z_{s-})}{dr} = \Delta L_{s}\,\kappa(\theta(r; \Delta L_{s}, z_{s-})), \quad \theta(0; \Delta L_{s}, z_{s-}) = z_{s-} \label{eqn:ivp_marcus}
\end{equation} 
\cite{Applebaum2004, Chechkin2014}. In the Marcus map $\theta(r; \Delta L_{s}, z_{s-})$,  $r$ is a time-like variable in which the process $z_{t}$ travels infinitely fast along a curve connecting the initial and end ``times'' of the jump, $\Delta L_{s}$, given by $r = 0$ and $r = 1$, respectively. The jump occurring at time $s$, denoted $\Delta L_{s} = dL_{s}^{(\alpha,\beta)}$, is a random variable distributed according to an $\alpha$-stable law, $\mathcal{S}_{\alpha}(\beta,dt^{1/\alpha})$. The resulting effect of the jump on the process $z_{t}$ is determined by the size of the jump as well as the integrated effect of the noise coefficient over the course of the jump \cite{Applebaum2004}. When an expression for $\theta$ can be derived and expressed in a closed form expression, then $z_{t}$ can be numerically estimated as
\begin{equation}
\hat{z}_{n+1} = \hat{z}_{n} + H(\hat{z}_{n})\,\delta t + [\theta(1;\delta L_{n}^{(\alpha,\beta)},\hat{z}_{n}) - \hat{z}_{n}]. \label{num_method_marcus_basic}
\end{equation}
where $\delta L_{n}^{(\alpha,\beta)}$ is as described in Appendix \ref{app:simulation_Ito}. For example, in the (N+) approximation of nonlinear system 1 (\ref{nonlin1_N+}), $\theta(1;\delta L_{n}^{(\alpha,\beta)},\hat{z}_{n}) = \hat{z}_{n}\exp\left(b\tau\,\delta L_{n}^{(\alpha,\beta)}\right)$.

In cases where $\theta$ needs to be evaluated numerically, as in the case of the (N+) approximation of nonlinear system 2 (\ref{nonlin2_N+}), we follow \cite{Asmussen2007book, Grigoriu2009} to approximate the $\alpha$-stable noise increments $dL_{t}^{(\alpha,\beta)}$ as the sum of a (path continuous) Gaussian white noise process and a (pure jump) compound Poisson process. This process effectively decomposes $dL_{s}$ into small jumps and large jumps, approximated respectively with the Gaussian white noise process and the Poisson process. We consider the case where the noise is symmetric (i.e. $\beta = 0$) as the case of asymmetric noise has not been studied explicitly (to our knowledge). It is necessary to define a large jump threshold $R > 0$ to approximate $dL_{t}^{(\alpha,0)}$ as
\begin{equation}
dL_{t}^{(\alpha,0)} \approx \eta\,dW_{t} + dN_{t} \label{eqn:stable_noise_approximation}
\end{equation}
where $\eta$ is a constant satisfying
\begin{equation}
\eta^{2} = \left(\frac{\alpha}{2 - \alpha}\right)C_{\alpha}R^{2 - \alpha}, \quad C_{\alpha} = \begin{cases} \frac{1 - \alpha}{\Gamma(2 - \alpha)\cos(\pi\alpha/2)} &\mbox{if $\alpha \ne 1$} \\ 2/\pi &\mbox{if $\alpha = 1$},\end{cases}
\end{equation}
$W_{t}$ is a standard Wiener process, and $N_{t} = N_{t}(\nu,\lambda)$ is a compound Poisson process with rate $\lambda$ and measure $\nu$ given by
\begin{align}
\lambda = C_{\alpha}/R^\alpha, \quad \nu(x) = \begin{cases}\frac{\alpha R^\alpha}{2}|x|^{-\alpha + 1} &\mbox{if $x < -R$ or $x \ge R$} \\
0 &\mbox{otherwise.}
\end{cases} \label{eqn:Poisson_params}
\end{align} 
We choose $R = 1$ to be our large jump threshold, as is done in \cite{Grigoriu2009}. Since we are approximating $\alpha$-stable noise with a sum of a continuous path process and pure jump process, we must use the general representation for stochastic integrals of L\'{e}vy processes that includes contributions from both the continuous and jump components \cite{Applebaum2004, Cont2004}, yielding
\begin{align}
\int_{0}^{t}\kappa(z_{s})\diamond dL_{s}^{(\alpha,0)} \approx &\,\,\eta\int_{0}^{t}\kappa(z_{s})\,dW_{s} + \frac{\eta^{2}}{2}\int_{0}^{t}\kappa(z_{s})\kappa'(z_{s})\,dt + \int_{0}^{t}\kappa(z_{s})\diamond dN_{s} \\ = &\,\,\eta\int_{0}^{t}\kappa(z_{s})\,dW_{s} + \frac{\eta^{2}}{2}\int_{0}^{t}\kappa(z_{s})\kappa'(z_{s})\,dt + \sum_{s \le t}\left[\theta(1; \Delta N_{s},z_{s-}) - z_{s-}\right] \label{marcus_discret_reduced}
\end{align}
where the first integral on the right-hand side is interpreted in the sense of It\={o} . The term $\Delta N_{s}$ represents the jump in the process $N_{t}$ occurring at time $s$ (if any), accounting for large jumps in the process $L_{t}^{(\alpha,0)}$. Then $\Delta N_{s}$ is distributed as a compound Poisson process evaluated over the infinitesimal time step $dt$ with rate and measure as in (\ref{eqn:Poisson_params}). We see that when $\alpha = 2$, then the jump rate $\lambda = 0$ (i.e. $\Delta N_{t} = 0$ for all $t$ and the sum in (\ref{marcus_discret_reduced}) vanishes). Thus (\ref{marcus_discret_reduced}) reduces to the sum of the It\={o} integral and the Wong-Zakai correction and $\eta^{2} = 2$, which is equivalent to the Stratonovich integral \cite{Oks03}. When $\alpha = 2$ and $\kappa$ is constant, then (\ref{marcus_discret_reduced}) reduces to the It\={o} integral, as expected. 

To simulate (\ref{eqn:Marcus_sde_app}) with $\beta = 0$, and time step $\delta t$, we use the numerical  scheme
\begin{equation}
\hat{z}_{n+1} = \hat{z}_{n} + \left(H(\hat{z}_{n}) + \frac{\eta^{2}}{2}\kappa(\hat{z}_{n})\kappa'(\hat{z}_{n})\right)\,\delta t + \eta \kappa(\hat{z}_{n})\,\delta W_{n} + [\hat{\theta}(1;\delta N_{n},\hat{z}_{n}) - \hat{z}_{n}] \label{num_method_marcus}
\end{equation}
where $\delta W_{n} \sim \mathcal{N}(0,\delta t)$ and $\delta N_{n}$ are discrete increments of a compound Poisson process over the time step $\delta t$ with rate and measure given by (\ref{eqn:Poisson_params}) which can be simulated as in \cite{Cont2004}. The term $\hat{\theta}(1;\delta N_{n},\hat{z}_{n}) = \hat{\theta}_{M}$ is the numerical approximation to $\theta(1;\delta N_{n},\hat{z}_{n})$, where
\begin{equation}
\hat{\theta}_{j + 1} = \hat{\theta}_{j} + \delta N_{n} \kappa(\hat{\theta}_{j})\,\delta u, \quad \hat{\theta}_{0} = \hat{z}_{n}, \quad \left(j = 0,1,\dots,M-1, \quad \delta u = \frac{1}{M} \right).  \label{eqn:theta_numerical}
\end{equation}
This numerical method for simulating $z_{t}$ when $\hat{\theta}$ needs to be evaluated is similar to the tau leaping methods that are given in \cite{Li2013} for evaluating Marcus integrals, and a more explicit description can be found therein (also see \cite{Li2013Erratum} and \cite{Sun2013}). Our numerical method converges weakly and more technical details related to the convergence of our method are discussed in \cite{Cont2004, Grigoriu2009}. 

Here we have written the numerical solutions (\ref{num_method_marcus_basic}), (\ref{num_method_marcus}), and (\ref{eqn:theta_numerical}) using the Euler method for the drift terms, but it would be possible to use higher-order methods, such as the explicit trapezoidal method to numerically integrate the drift terms.

\subsection{Estimating the PDFs}
\label{app:stats}

To ensure accurate estimation of the tails of the sample PDFs, we estimate a value $N$ for the size of the sample of the time varying processes that provides reasonable confidence in the estimates for the PDFs. This value of $N$ is based on a calculation for estimated confidence intervals for the tails of the PDFs. Here $N$ is equal to the length of the time interval of the simulation divided by $Dt$, so that in our calculations we have $N$ discrete observations of process $z_t$ sampled at intervals of $Dt=10^{-2}$.  However, in order to obtain an approximation for the confidence intervals of the tails of the PDF, we use approximations based on independent observations, so we require a subsample $Z_j,\, j = 1,\dots, N_Z$ of $z_t$ which has $N_Z$ almost-independent observations. Given the sampling time step $Dt  = 10^{-2}$ and that the characteristic time scale of the slow variable is $O(1)$ (as can be seen by inspection of the AFs above in Fig. \ref{lin_autocodiff}, \ref{nonlin1_autocodiff}), we need to subsample $z_{t}$ once every $1/Dt = 10^{2}$ points to ensure that observations of $Z$ are sufficiently independent.  Then $N$ is determined as $N_Z \times 10^2$, with $N_Z$ determined from confidence intervals for the PDF of $Z$.

To obtain confidence intervals for the tails of the PDFs, we discretize the state space of $Z$ into $n$ disjoint bins where the true probability of observing $Z_{j}$ within bin $i$ is $\pi_{i}$. Under these assumptions, the estimated number of observations of $Z_{j}$ in the $n$ discrete bins is given by a multinomial distribution with $N_{Z}$ observations and probabilities $\{\pi_{i}\}_{i = 1}^{n}$ ($\sum_{i}^{n} \pi_{i} = 1$). Obviously, the $\pi_{i}$ for the bins that correspond to the tail of the density are very small, yet we want to ensure that $N_{Z}$ (and thus $N$) is sufficiently large to obtain a good approximation for the tails of the density in our simulations. Therefore, we derive an expression for the confidence intervals of $\pi_i$ based on $N_Z$ and $n$. Specifying that small $\pi_{i}$ ($\pi_{i} = 10^{-5}$) fall within a 95\% confidence interval of width $O(10^{-5})$ about $\pi_i$, we obtain a minimum value of $N_{Z}$. First, we approximate the multinomial distribution for the number of observations using a multivariate normal approximation which is justified on the basis of large $N_{Z}$. We then apply the confidence interval suggested in \cite{Miller1981} by Quesenberry and Hurst with Goodman's recommendation and apply the limit $N_{Z} \rightarrow \infty$ to obtain the following approximate two-sided $95\%$ confidence interval for $\pi_{i}$,
\begin{align}
&\prob{\hat{\pi}_{i} - \Omega_{i} < \pi_{i} < \hat{\pi}_{i} + \Omega_{i}} \approx 0.95 \label{QH_interval_asymp} \\
&\hat{\pi}_{i} = \frac{1}{N}\left\{\mbox{\# observations of $Z \in$ bin $i$} \right\}, \\ 
&\Omega_{i} = \sqrt{\frac{M_{Z}\hat{\pi}_{i}(1 - \hat{\pi}_{i})}{N_{Z}}}, \quad M_{Z} = \chi^{2}_{1,1 - 0.05/n}
\end{align}
where $\chi^{2}_{m,1-c}$ denotes the $(1-c)\%$ quantile of the $\chi^{2}$ distribution with $m$ degrees of freedom. To get a conservative estimate of the value of $\Omega_{i}$ we might expect for our simulations, we choose $n = 100$, and so $M_{Z} = \chi^{2}_{1,1 - 0.05/100} \approx 12.12 = O(10^{1})$. We wish to have narrow confidence intervals where $\hat{\pi}_{i}$ is within $O(10^{-5})$ of $\pi_{i} = 10^{-5}$, which implies that $\Omega_{i} = O(10^{-5})$ and $N_{Z}$ satisfies,
\begin{equation}
N_{Z} = M_{Z}\pi_{i}/\Omega_{i}^{2} = O(10^{6}). \label{eqn:Nz}
\end{equation}
Since $N = N_{Z}/Dt$, we conclude that $10^{8}$ points with time step $Dt$ are needed in the full simulation to resolve the characteristic time scale, the sampling of $z_{t}$, and the subsampling of $Z$.

\subsection{Estimating the autocodifference function}
\label{app:numericalACD}
Given a time-dependent stochastic process $z_{t}$, the autocodifference of $z_{t}$, $A_{z}(\tau)$, is defined by (\ref{eqn:autocodiff}) \cite{Rosadi2009, Taqqu1994}. To numerically estimate $A_{z}(s)$ we need to estimate both the CF of $z_{t}$, $\Phi_{z}(k) = \expected{\exp(ikz_{t})}$ and the joint CF of $z_{t}$ and a time-shifted version, $z_{t + \tau}$, $\Psi_{z}(k,l,\tau) = \expected{\exp(ikz_{t+\tau} +ilz_{t})}$. Assume we have $N$ observations of $z_{t}$, which we denote by $\hat{z}_{n} = z_{n\cdot Dt}$, $n = 0,1,\dots, N-1$. The estimates for $\Phi_{z}(k),\,\Psi_{z}(k,l,n\cdot Dt)$ are given by
\begin{equation}
\begin{cases}
\hat{\Phi}_{z}^{[J,K]}(k) = \displaystyle \frac{1}{K - J + 1}\left(\sum_{j = J}^{K}\exp(ik\hat{z}_{j})\right), \\ 
\hat{\Psi}_{z}^{[J,K]}(k,l,n) = \displaystyle \frac{1}{K - J + 1}\left(\sum_{j = J}^{K}\exp(ik\hat{z}_{j + n} + il\hat{z}_{j})\right),
\end{cases}
\label{eqn:char_func_ests}
\end{equation}
respectively, where $J,K$ are integers satisfying $0 \le J \le K \le N-1-n$. The estimate for $A_{z}(n\cdot Dt)$ which we denote by $\hat{A}_{z}(n\cdot Dt)$ can be computed at discrete time displacements $n\cdot Dt$ using (\ref{eqn:char_func_ests}) and is given by
\begin{equation}
\hat{A}_{z}(n\cdot Dt) = \log\left(\hat{\Psi}_{z}^{[0,N-1-n]}(1,-1,n)\right) - \log\left(\hat{\Phi}_{z}^{[n,N-1]}(1)\right) - \log\left(\hat{\Phi}_{z}^{[0,N-1-n]}(-1)\right)
\end{equation}
This estimator is reasonably accurate for short times where $n \ll N$ as there are $N - n$ pairs of data points with relative time-displacement $n\cdot \delta t$, and by the law of large numbers $\hat{\Psi}$, $\hat{\Phi}$ given in (\ref{eqn:char_func_ests}) are reasonable estimators. For our study, we consider $n \leq 400 \ll N$ and so $A(n\cdot Dt)$ is well-estimated.

\section{Acknowledgements}
The authors wish to acknowledge the generous support of NSERC\footnote{National Science and Engineering Research Council (Canada): \url{http://www.nserc-crsng.gc.ca}}. The authors thank Dr. I. Pavlyukevich (Uni. Jena) for bringing the Marcus integral to our attention and for helpful advice. WFT also thanks Dr. J. Walsh (UBC)  for helpful discussions.

\bibliographystyle{siam}
\bibliography{phd}

\end{document}